\documentclass[12pt,draft,leqno]{article}
\usepackage{amssymb,eucal,latexsym}
\usepackage{amsmath,amsthm,latexsym,amscd,amsbsy,fleqn,graphicx,upgreek}

\newcommand{\dd}{\partial}

\newcommand{\f}{\ensuremath\varphi}
\newcommand{\F}{\ensuremath\mathcal{F}}
\newcommand{\G}{\ensuremath\mathcal{G}}
\newcommand{\PP}{\ensuremath\mathcal{P}}

\newcommand{\Rn}{\mathcal{R}^n}

\newcommand{\la}{\langle }
\newcommand{\ra}{\rangle }
\newcommand{\lh}{\leftrightharpoons}
\newcommand{\RRnm}{$\mathbb{R}^n~$}

\newcommand{\Cl}{\textit{Cl}}
\newcommand{\Int}{\textit{Int}}
\newcommand{\Ra}{\textit{Range}}
\newcommand*{\twoheadrightarrowtail}{\mathrel{\rightarrowtail\kern-1.9ex\twoheadrightarrow}}

\newtheorem{theorem}{Theorem}[section]
\newtheorem{lm}[theorem]{Lemma}
\newtheorem{exa}[theorem]{Example}

\newtheorem{cor}[theorem]{Corollary}

\newtheorem{defi}[theorem]{Definition}

\newtheorem{rem}[theorem]{Remark}

\title{{\LARGE\bf A Logic of Strong Contact}\\
\vspace{0.2cm}
{\LARGE\bf  between Polytopes}\\
\vspace{0.5cm}
{\large{\bf Tsvetlin Marinov}
 and {\bf Tinko Tinchev}}\\
\vspace{0.2cm}
{\footnotesize\rm Department of Mathematical Logic and Applications, Sofia University St. Kliment Ohridski}\\
{\footnotesize\rm 5 J. Bourchier Blvd., 1164 Sofia, Bulgaria}
}

\author{}

\date{}

\begin{document}

\maketitle

\begin{abstract}
We propose a new contact relation between polytopes. Intuitively, we say that two polytopes are in strong contact if a small enough object can pass from one of them to the other while remaining in their union. In the first half of the paper we prove that this relation is indeed a contact relation between polytopes, which turns out not to be the case for arbitrary regular closed in Euclidean spaces sets. In the second half we study the universal fragments of the logics of the resultant contact algebras. We prove that they all coincide with the set of theorems of a standard quantifier-free formal system for connected contact algebras, which also coincides with the universal fragments of the logics of a variety of (classes of) contact algebras of interest.
\end{abstract}

\footnotetext[1]{{\footnotesize
{\em Keywords:}  polytope, contact relation, contact algebra, adjacency space, logic of space}}

\footnotetext[2]{{\footnotesize {\em E-mail addresses:}
cvmarinov@fmi.uni-sofia.bg, tinko@fmi.uni-sofia.bg}}

\section{Introduction}

\textit{Region-based theory of space} (RBTS) is an alternative to the standard \textit{point-based theory of space}. It originates from the philosophical argument, proposed by Whitehead \cite{whitehead1957process}, de Laguna \cite{deLaguna1922point} and others, that the notion of a \textit{point} is too abstract to be taken as primitive. They reasoned that the primitive ontological notion of geometry should, instead, resemble \textit{spatial bodies}, for which the name \textit{region} has been chosen, and that the notion of a point should be \textit{defined} in terms of the notion of region and basic relational notions such as \textit{part-of} and \textit{contact}. In addition to this philosophical motivation, this approach to geometry has been a subject of interest due to its applicability in practical areas such as qualitative spatial reasoning (QSR), knowledge representation, geographical information systems, formal ontologies in information systems, image processing and natural language semantics.\\

Early papers related to this idea include Tarski's \cite{tarski1929fondements,tarski1956foundations} , Grzegorczyk's \cite{grzegorczyk1961axiomatizability}, Clarke's \cite{clarke1981calculus} and \cite{clark1985individuals} and Biacino and Gerla's \cite{biacino1991connection}. Some more recent works, focused on the correspondence with the point-based approach, are Roeper's \cite{roeper1997region}, Mormann's \cite{mormann1998continuous}, Pratt and Schoop's \cite{pratt1998complete}, Vakarelov, Dimov, D\"{u}ntsch and Bennet's \cite{vakarelov2001note,vakarelov2002proximity}, D\"{u}ntsch and Winter's \cite{duntsch2005representation} and Dimov and Vakarelov's \cite{dimov2006contact,dimov2006contact2}. An extensive survey in Spatial Logic is \cite{aiello2007handbook}. One of the important systems of RBTS is the Region Connection Calculus (RCC), introduced by Randel, Cui and Cohn in \cite{randell1992spatial}, for which an axiomatisation based on Boolean algebras was given by Stell in \cite{stell2000boolean}. It and a number of modifications of RCC have been intensively studied \cite{egenhofer1991point,smith1992algebraic,kontchakov2014spatial,kontchakov2013topological}.\\

Different objects could be taken as regions. A standard choice is the regular closed sets of suitable topological spaces, like Euclidean spaces. But such sets can have very exotic properties, like, for instance, some fractals. One possible restriction to more \textit{tame} sets, which presumably better resemble spatial bodies, are the \textit{polytopes}, which we consider, which are a special kind of regular closed in Euclidean spaces sets.\\

Also, different contact relations can be considered. In this paper we propose a new kind of contact relation between polytopes, which seems to have good resemblance with our natural perception of the notion of contact. We say, intuitively, that two polytopes are in strong contact if a sufficiently small object can pass from one of them to the other while remaining in their union. This idea has a nice topological formalisation. The resultant relation is strictly stronger than the standard topological contact (for all dimensions greater than 1) and strictly weaker than the overlap relation.\\

In addition to the short preliminary Section 1, this paper has two main parts. In Section 2 we define the strong-contact relation and prove that it is indeed a contact relation in the considered algebras of polytopes. Interestingly, it turns out that it is not a contact relation between arbitrary regular closed in Euclidean spaces sets. In Section 3 we study the universal fragments of the logics of the resultant \textit{strong-contact polytope algebras}.
%To do that we consider a standard quantifier-free formal system for connected contact algebras and show that its theorems are precisely the true formulas in each of these algebras. 
To do that we prove completeness theorems for a standard quantifier-free formal system for connected contact algebras with respect to particular structures of interest.\\

\section{Preliminaries}

\subsection{Boolean Algebras}

Let $A$ be a nonempty set, $-$ be a unary operation in $A$, $+$ and $\cdot$ be binary operations in $A$ and $0$ and $1$ be two distinct elements of $A$. Let for any elements $x$, $y$ and $z$ of $A$ the following conditions be satisfied:

\begin{center}
	\begin{tabular}{ r c l }
		$x+(y+z)=(x+y)+z$ 					& (associativity) 	& $x\cdot(y\cdot z)=(x\cdot y)\cdot z$\\
		$x+y=y+x$ 							& (commutativity) 	& $x\cdot y=y\cdot x$\\
		$x+(x\cdot y)=x$ 					& (absorption) 		& $x\cdot(x+y)=x$\\
		$x\cdot(y+z)=(x\cdot y)+(x\cdot z)$ & (distributiviy) 	& $x+(y\cdot z)=(x+ y)\cdot(x+z)$\\
		$x+(-x)=1$ 							& (complementation) & $x\cdot(-x)=0$\\
	\end{tabular}
\end{center}

\noindent Then $\mathcal{A}=\la A,-,+,\cdot,0,1\ra$ is called a \textit{Boolean algebra} and $-$, $+$, $\cdot$, $0$ and $1$ are called respectively the \textit{complement}, \textit{join}, \textit{meet}, \textit{bottom element} (or \textit{zero}), and \textit{top element} (or \textit{unit}) of $\mathcal{A}$. The binary relation $\leq$ in $A$, such that $x\leq y$ iff $x+y=y$, is called \textit{the Boolean ordering of $\mathcal{A}$}.\\

By the complementation equalities and the known fact that the complement, join and meet satisfy de Morgan's laws, any Boolean algebra is determined by its carrier, complement and join. That is why when we say "the Boolean algebra $\la A,-,+\ra$" we mean the unique Boolean algebra with carrier $A$, complement $-$ and join $+$.\\

Let $W$ be a nonempty set. Then the power set $\mathcal{P}(W)$ of $W$ is the carrier of a Boolean algebra with complement the set-theoretic complement $W\setminus~$ to $W$ and join the set-theoretic union $\cup$. In other words $\la\PP(W),W\setminus~,\cup\ra$ is a Boolean algebra. We shall designate it by $B(W)$ and call it the \textit{set-theoretic Boolean algebra over $W$}. Its meet, zero, unit and ordering are respectively the set-theoretic intersection $\cap$, the empty set $\emptyset$, the set $W$ and the set-theoretic inclusion $\subseteq$.\\

Let $\mathcal{A}=\la A,-,+\ra$ be a Boolean algebra. If $B$ is a closed with respect to $-$ and $+$ nonempty subset of $A$, we say that $\mathcal{B}=\la B,-,+\ra$ is a \textit{subalgebra of $\mathcal{A}$}. Clearly a subalgebra $\mathcal{B}$ of a Boolean algebra $\mathcal{A}$ is itself a Boolean algebra and its meet, zero and unit are respectively the meet, zero and unit of $\mathcal{A}$.\\

The notions of join and meet are generalised to arbitrary nonempty subsets of the carrier of a Boolean algebra. Let $\mathcal{A}$ be a Boolean algebra with carrier $A$ and ordering $\leq$. Let $B$ be a nonempty subset of $A$. An element $+B$ of $A$ is said to be \textit{the join in $\mathcal{A}$ of $B$} iff $(\forall b\in B)(b\leq +B)$ and $(\forall a\in A)((\forall b\in B)(b\leq a)~\rightarrow~ +\!B\leq a)$. Analogically, an element $\cdot B$ of $A$ is said to be \textit{the meet in $\mathcal{A}$ of $B$} iff $(\forall b\in B)(\cdot B\leq b)$ and $(\forall a\in A)((\forall b\in B)(a\leq b)~\rightarrow~ a\leq\cdot B)$. A Boolean algebra is said to be \textit{complete} iff its carrier contains the join and meet of each of its nonempty subsets.\\

\subsection{Contact Relations and Contact Algebras}

Let $\mathcal{A}$ be a Boolean algebra with carrier, complement, join, meet, zero, unit and ordering respectively $A$, $-$, $+$, $\cdot$, $0$, $1$ and $\leq$.

\begin{defi}\label{contact relation}\rm
	A binary relation $C$ in $A$ is called a \textit{contact relation in $\mathcal{A}$} iff, for any elements $x$, $y$ and $z$ of $A$, the following conditions are satisfied:
	\begin{center}
		\begin{tabular}{ c l }
			(C1) ~&~ $ \neg C(0,x) $\\
			(C2) ~&~ $ C(x,y+z)\leftrightarrow (C(x,y) ~\textit{or}~ C(x,z))$\\
			(C3) ~&~ $ C(x,y)\rightarrow C(y,x)$\\
			(C4) ~&~ $ x\neq 0\rightarrow C(x,x)$
		\end{tabular}
	\end{center}
\end{defi}

\begin{defi}\label{contact algebra}\rm
	If $\mathcal{A}$ is a Boolean algebra and $C$ is a contact relation in $\mathcal{A}$, then $\la \mathcal{A},C\ra$ is called a \textit{contact algebra}.
\end{defi}

\begin{exa}\rm
	We say that two elements $x$ and $y$ of $A$ \textit{overlap} iff their meet is not the bottom element, i.e. $x\cdot y\neq 0$. It is easy to see that the overlap relation in any Boolean algebra is a contact relation in it.
\end{exa}

\begin{defi}\label{adjacency space}\rm

Let $W$ be a nonempty set and $R$ be a binary relation in $W$. Let $C_R$ be the binary relation in $\PP(W)$ such that for any subsets $a$ and $b$ of $W$ we have $C_R(a,b)$ iff $(\exists x\in a)(\exists y\in b)xRy$. It is easy to see that $\la B(W),C_R\ra$ is a contact algebra iff $R$ is reflexive and symmetric. If that is the case, we call $\F=\la W,R\ra$ an \textit{adjacency space}, we call the elements of $W$ \textit{cells of $\F$}, we call $R$ the \textit{adjacency relation of $\F$} and we say that the contact algebra $\la B(W),C_R\ra$ is \textit{induced by $\F$}. We call a contact algebra which is induced by some adjacency space a \textit{set-theoretic contact algebra}.

\end{defi}

The following properties of contact relations are well-known and follow easily from the conditions (C1) to (C4). \\
A contact relation is monotone with respect to the Boolean ordering, i.e. $$x\leq x'\rightarrow(y\leq y'\rightarrow(C(x,y)\rightarrow C(x',y')))~.$$ Any contact relation is an extension of the overlap relation, i.e. if two elements overlap, they are in contact, i.e. $$x\cdot y\neq0\rightarrow C(x,y)~.$$\\

\subsection{Topological Contact}

Let $T=\la X,\tau\ra$ be an arbitrary topological space.\\

Let $\Int$, $\Cl$ and $\partial$ designate respectively the interior, closure and boundary operators. Let $\sqcap$ designate the binary operation, called \textit{regularised intersection}, such that for any subsets $A$ and $B$ of $X$ we have $A\sqcap B\lh \Cl(\Int(A\cap B))$. Let $*$ designate the unary operation such that for any subset $A$ of $X$ we have $A^*\lh\Cl(X\setminus A)$.\\%is the closure of the set-theoretic complement of $A$ to the universe $X$ of $T$.\\

A subset $A$ of $X$ is called \textit{regular closed in $T$} iff $A=\Cl(\Int(A))$. We designate the set $\{A\subseteq X\mid A=\Cl(\Int(A))\}$ of the regular closed in $T$ sets by $RC(T)$. It is known, for instance from \cite{sikorski1969boolean}, that $\mathcal{RC}(T)\lh\la RC(T),*,\cup\ra$ is a complete Boolean algebra with meet, zero, unit and ordering respectively $\sqcap$, $\emptyset$, $X$ and $\subseteq$. Moreover, the join and meet of a set $A$ of regular closed sets equal $\Cl(\cup A)$ and $\Cl(\Int(\cap A))$ respectively. In particular, if $A=\{a_1,...,a_k\}$ is a finite set of regular closed sets, we have $a_1\sqcap...\sqcap a_k=((...((a_1\sqcap a_2)\sqcap a_3)\sqcap...)\sqcap a_{k-1})\sqcap a_k=\Cl(\Int(\cap A))=((...((a_1\cap a_2)\cap a_3)\cap...)\cap a_{k-1})\sqcap a_k$, which we shall designate by $\sqcap A$ and call the \textit{regularised intersection of $A$}. \\

Let us point out that, since for any set $B$ in a topological space we have $\dd\Cl(B)\subseteq\dd B$, we have that the boundary points of a regular closed set $A$ are boundary points of its interior and thus any open neighbourhood of such a point contains not only points of $A$ but points of $\Int(A)$ as well.\\

Let $A$ and $B$ be regular closed in $T$ sets. We say that $A$ and $B$ are in \textit{topological contact} iff $A\cap B\neq\emptyset$. We shall designate this binary relation by $C^T$, or just $C$ for short. It is easy to verify that $C^T$ is a contact relation in $\mathcal{RC}(T)$.\\

\section{Strong Contact}

Let $T=\la X,\tau\ra$ be a topological space. We say that an open in $T$ set is \textit{connected} iff it cannot be represented as the union of two disjoint open sets. Let us define the binary relation $SC^T$ in $\mathcal{P}(X)$ as follows: 

\begin{defi}\label{strong contact}\rm
	For any subsets $A$ and $B$ of $X$, let $SC^T(A,B)$ iff there exists a connected and open subset $E$ of $A\cup B$ such that $E\cap A\neq\emptyset$ and $E\cap B\neq\emptyset$.
\end{defi}

 We shall omit the superscript when it is clear from the context.

\begin{rem}\label{ContIndeed}\rm
	Obviously $SC$ is symmetric and $\neg SC(\emptyset,A)$ for any $A$. Also, evidently $SC(A,B)$ or $SC(A,D)$ implies $SC(A,B\cup D)$.
\end{rem}

We shall consider the $SC$ relations for Euclidean spaces. Let for any positive natural number $n$, $\mathcal{R}^n$ be the set of $n$-tuples of real numbers, $\mathcal{T}^n$ be the natural topology on $\mathcal{R}^n$ and $\mathbb{R}^n=\la\mathcal{R}^n,\mathcal{T}^n\ra$.\\

\begin{lm}\label{upward strength}\textsc{(Strength, upward.)}
	$SC^{\mathbb{R}^n}$ is an extension of the overlap relation in $\mathcal{RC}(\mathbb{R}^n)$.
\end{lm}

\textit{Proof}. Let $A$ and $B$ be overlapping regular closed in $\mathbb{R}^n$ sets. Then $A\sqcap B\neq\emptyset$, i.e. $\Cl(\Int(A\cap B))\neq\emptyset$, thus $\Int(A\cap B)\neq\emptyset$. Let $x\in\Int(A\cap B)$. Evidently, any open ball with centre $x$, contained in $\Int(A\cap B)$ is a witness to $SC(A,B)$.

\begin{cor}\rm
	\label{UpwardStrengthCorollary}
	Evidently, for any nonempty regular closed set $A$, we have $SC(A,A)$.\\
\end{cor}

By this and remark \ref{ContIndeed}, to prove that $SC$ is a contact relation in a subalgebra of $\mathcal{RC}(\mathbb{R}^n)$, it remains only to prove the left-to-right direction of (C2), i.e. that $SC$ is distributive over the join of the algebra in question, i.e. that $SC(A,B\cup D)$ implies $SC(A,B) ~or~ SC(A,D)$.\\

\begin{lm}\label{downward strength}\textsc{(Strength, downward.)}
	Let $A$ and $B$ be closed in $\mathbb{R}^n$ sets such that $SC(A,B)$. Then $A\cap B\neq\emptyset$.
\end{lm}

\textit{Proof}. Let $E$ be a witness to $SC(A,B)$. Suppose $A\cap B=\emptyset$. Since $\mathbb{R}^n$ is a normal topological space, let $A'$ and $B'$ be open sets such that $A\subseteq A'$, $B\subseteq B'$ and $A'\cap B'=\emptyset$. Then $E\subseteq A\cup B\subseteq A'\cup B'$, so $E=(E\cap A')\cup(E\cap B')$. Thus $E$ is the union of two nonempty disjoint open sets, i.e. $E$ is not connected, which is a contradiction.\qed\\

We shall show that the relation $SC^{\mathbb{R}^n}$ is not distributive over the set-theoretic union for regular closed sets, by showing a counterexample in $\mathcal{RC}(\mathbb{R}^1)$. We shall use the partitioning of the closed interval $[0,1]$ by the sequence of the negative integer powers of 2. Let for any natural number $k$, $S_k$ designate the closed interval $[2^{-k-1},2^{-k}]$.\\

Let $B\lh\Cl(\cup\{S_{2k}\mid k<\omega\})$ be the closure of the union of those line segments $S_k$ with even indices and $D\lh\Cl(\cup\{S_{2k+1}\mid k<\omega\})$ -- of those with odd indices. $B$ and $D$ are defined as the joins in $\mathcal{RC}(\mathbb{R}^1)$ of $\{S_{2k}\mid k<\omega\}$ and $\{S_{2k+1}\mid k<\omega\}$ thus $B$ and $D$ are regular closed in $\mathbb{R}^1$ sets.\\

The point 0 is an accumulation point of both $B$ and $D$, thus $0\in\Cl(B)=B$ and $0\in\Cl(D)=D$. Clearly $B\cup D=[0,1]$. Let also $A\lh[-1,0]$.\\

The open interval $(-1,1)$ is a witness to $SC(A,B\cup D)$. But $\neg SC(A,B)$ because no open interval (connected and open in $\mathbb{R}^1$ set) which has nonempty intersection with both $A$ and $B$ is contained in $A\cup B$. Analogically $\neg SC(A,D)$. Thus $SC$ is not distributive over $\cup$ for regular closed sets in $\mathbb{R}^1$, thus $SC^{\mathbb{R}^1}$ is not a contact relation in $\mathcal{RC}(\mathbb{R}^1)$.\\

\subsection{Polytopes}

We shall now define a particular kind of regular closed in Euclidean spaces sets, which we shall call polytopes.

\begin{defi}[]\label{polytope}\rm
	A regularised intersection of finitely many closed half-spaces of $\mathbb{R}^n$ is called a \textit{basic polytope in $\mathbb{R}^n$}. 	A finite union of basic polytopes in $\mathbb{R}^n$ is called a \textit{polytope in $\mathbb{R}^n$}. 
\end{defi}

We shall designate the set of polytopes in $\mathbb{R}^n$ by $P^n$.

\begin{rem}\rm
	Notice that $\emptyset$ and $\mathcal{R}^n$ are polytopes in $\mathbb{R}^n$, since for any closed half-space $\alpha$ we have that $\alpha^*$ is also a closed half-space and $\alpha\sqcap\alpha^*=\emptyset$ and $\alpha\cup\alpha^*=\mathcal{R}^n$. Notice also that polytopes are regular closed.
\end{rem}

A set $A$ in an Euclidean space is called \textit{convex}, iff each line segments with endpoints belonging to $A$ is a subset of $A$. \\

We shall use the following well-known, described, for instance, in \cite{de2000mathematical}, results about convex sets in Euclidean spaces: a closed half-space of a Euclidean space is convex; the intersection of any set of convex sets is convex; if $A$ is a convex set with nonempty interior, then $\Cl(A)=\Cl(\Int(A))$. Using these results we immediately obtain the following

\begin{lm}\label{basic polytopes}\textsc{(Basic polytopes.)}
	If $A=\sqcap B$ is a nonempty basic polytope, then $\cap B=\Cl(\cap B)=\Cl(\Int(B))=\sqcap B$.
\end{lm}

We shall now show that the polytopes in $\mathbb{R}^n$ form a Boolean subalgebra of $\mathcal{RC}(\mathbb{R}^n)$, i.e. that $P^n$ is closed with respect to the operations $*$ and $\cup$. Clearly the union of two polytopes is a polytope.\\

Let $A$ be a polytope. We shall show that $A^*$ is also a polytope. Let $A=\cup_{i=1}^qA_i=\cup_{i=1}^q(\sqcap_{j=1}^{p_i}\alpha_{ij})$, where all $\alpha_{ij}$ are closed half-spaces and thus all $A_i$ are basic polytopes. By de Morgan's laws we have that $A^*=(\cup_{i=1}^qA_i)^*= \sqcap_{i=1}^q(A_i^*)=\sqcap_{i=1}^q((\sqcap_{j=1}^{p_i}\alpha_{ij})^*)=\sqcap_{i=1}^q(\cup_{j=1}^{p_i}(\alpha_{ij}^*))$.
So $A^*$ is a finite regularised intersection of finite unions of closed half-spaces. Thus it is a finite regularised intersection of finite unions of basic polytopes, i.e. a finite regularised intersection of polytopes. We shall prove that a regularised intersection of any finite number $q$ of polytopes is a polytope by induction on $q$.\\

Let every regularised intersection of $q$ polytopes be a polytope. Let $A^*=B\sqcap D_1\sqcap D_2\sqcap...\sqcap D_q$ be a regularised intersection of $q+1$ polytopes. If $q=0$, then $A^*=B$ is obviously a polytope, so let $q>0$. Since $B$ and $D_1$ are polytopes, let $B=B_1\cup...\cup B_s$ and $D_1=G_1\cup...\cup G_t$ for some basic polytopes $B_1$,..., $B_s$, $G_1$,...,$G_t$. Using simple properties of Boolean operations we obtain $B\sqcap D_1=B\sqcap(G_1\cup...\cup G_t)=(B\sqcap G_1)\cup...\cup (B\sqcap G_t)=((B_1\cup...\cup B_s)\sqcap G_1)\cup...\cup ((B_1\cup...\cup B_s)\sqcap G_t)=
(B_1\sqcap G_1)\cup...\cup(B_1\sqcap G_t)\cup...\cup(B_s\sqcap G_1)\cup...\cup(B_s\sqcap G_t)$. Thus $B\sqcap D_1$ is a finite union of regularised intersections of basic polytopes, thus is a finite union of basic polytopes, thus is a polytope. Then $A^*=B\sqcap D_1\sqcap D_2\sqcap...\sqcap D_q=(B\sqcap D_1)\sqcap D_2\sqcap...\sqcap D_q$, is a regularised intersection of $q$ polytopes and by the induction hypothesis is a polytope.\\

We shall designate the Boolean algebra $\la P^n,*,\cup\ra$ by $\mathcal{P}^n$.\\

\subsection{The One-dimensional Case}

By the downward strength lemma \ref{downward strength} we have that for any closed in $\mathbb{R}^1$ sets $A$ and $B$ we have that $SC^{\mathbb{R}^1}(A,B)$ implies $C^{\mathbb{R}^1}(A,B)$. We shall now show that if $A$ and $B$ are polytopes in $\mathbb{R}^1$ we also have that $C^{\mathbb{R}^1}(A,B)$ implies $SC^{\mathbb{R}^1}(A,B)$.\\

Polytopes in $\mathbb{R}^1$ are finite unions of closed intervals (with nonzero length) and/or rays. Let $A$ and $B$ be such and let $x\in A\cap B$. Let $A'$ and $B'$ be closed intervals (each with nonzero length) contained in $A$ and $B$ respectively such that $x$ is an endpoint of both. Let $a$ and $b$ be their other endpoints. Without loss of generality (WLoG) let $a\leq b$. If $x$ is between $a$ and $b$, then the open interval $(a,b)$ is a witness to $SC(A,B)$. If $ a\leq b<x$, then the open interval $(b,x)$ is a witness to $SC(A,B)$.\\

Thus for polytopes in $\mathbb{R}^1$ the relations $C^{\mathbb{R}^1}$ and $SC^{\mathbb{R}^1}$ coincide. Thus $SC^{\mathbb{R}^1}$ is a contact relation in $\PP^1$.\\

\subsection{The Two-dimensional Case}

It is easy to see that there are polytopes in $\mathbb{R}^n$ which have nonempty intersection but are not in the $SC^{\mathbb{R}^n}$ relation -- for instance a pair of vertical (opposite) angles in $\mathbb{R}^2$.\\

\begin{lm}\label{crossing}\textsc{(Crossing.)}
	Let $T=\la X,\tau\ra$ be a topological space, $A$ be a closed in $T$ set, $a\in\Int(A)$, $b\notin A$ and $\gamma$ be a curve in $T$ connecting $a$ and $b$ ($\gamma:[0,1]\longrightarrow X$, $\gamma(0)=a$ and $\gamma(1)=b$). Then $\Ra(\gamma)\cap\dd A\neq\emptyset$.
\end{lm}

\textit{Proof}. Let $B\lh\Cl(X\setminus A)=(X\setminus \Int(A))$. We will recursively define a sequence $\{x_i\}_{i<\omega}$ of points on $[0,1]$, as follows:\\
Base: $x_0\leftrightharpoons 0$\\
%Recursion hypothesis: Let $x_i$ be defined.\\
Recursion step:\\
If $\gamma(x_i)\in A\setminus B$, then let $x_{i+1}\leftrightharpoons x_i+2^{-i}$.\\
If $\gamma(x_i)\in B\setminus A$, then let $x_{i+1}\leftrightharpoons x_i-2^{-i}$.\\
If $\gamma(x_i)\in A\cap B$, then let $x_{i+1}\leftrightharpoons x_i$.\\

Notice that, since $x_0=0$, we have $\gamma(x_0)\in A\setminus B$, so $x_1=1$, so there exists $i$ such that $\gamma(x_i)\in B\setminus A$.

Case 1: $\exists i(x_{i+1}=x_i)$. Let $k$ be such. Then $\gamma(x_k)\in A\cap B=\dd A$.

Case 2: $\forall i(x_{i+1}\neq x_i)$. Then $\neg\exists i(\gamma(x_i)\in A\cap B)$

Suppose that only finitely many elements of $\{\gamma(x_i)\mid i<\omega\}$ belong to $A\!\setminus B$ and let $\gamma(x_k)$ be the last such (i.e. the one with the greatest index). Then $(\forall j>k)(\gamma(x_j)\in B\setminus A)$. Then $x_{k+1}=x_k-2^{-k}$ and for each $i>k$, we have $x_{i+1}=x_i+2^{-i}$. Then $$\lim_{i\rightarrow\omega}x_i~=~
x_k+\bigg(\frac{1}{2}\bigg)^{k} - \bigg(\frac{1}{2}\bigg)^{k+1} - \bigg(\frac{1}{2}\bigg)^{k+2}-...~=~
x_i+\bigg(\frac{1}{2}\bigg)^{k}-\sum_{i=1}^{\omega}\bigg(\frac{1}{2}\bigg)^{k+i}~=~x_k$$ By the continuity of $\gamma$, every open neighbourhood of $\gamma(x_k)$ contains a point $\gamma(x_{k+i})$ of $B$, thus $\gamma(x_k)\in \Cl(B)=B$, which contradicts $\gamma(x_k)\in A\!\setminus B$. Thus, infinitely many elements of $\{\gamma(x_i)\mid i<\omega\}$ belong to $A\!\setminus B$. Analogically, infinitely many elements of $\{\gamma(x_i)\mid i<\omega\}$ belong to $B\!\setminus A$. 

Since any series $\sum_{i<\omega}(-1)^{\epsilon(i)}2^{-i}$, where $\epsilon:\omega\longrightarrow\{0,1\}$, of the powers of $\frac{1}{2}$ is absolutely convergent, the sequence $\{x_i\}_{i<\omega}$ converges. Let $x\lh\lim_{i\rightarrow\omega}x_i$. Since $\gamma$ is continuous, we have that $\lim_{i\rightarrow\omega}\gamma(x_i)=\gamma(\lim_{i\rightarrow\omega}x_i)=\gamma(x)$. 
Let $\{a_i\}_i$ and $\{b_i\}_i$ be the subsequences of $\{x_i\}_i$ of those $x_i$ which are elements of $A\setminus B$ and those which are elements of $B\setminus A$ respectively. Then $x_k$ is a point of accumulation of both of them. Then every open neighbourhood of $\gamma(x_k)$ contains points form $A$ and points from $B$. Thus $\gamma(x_k)\in\Cl(A)\cap\Cl(B)=%\Cl(A)\setminus\Int(A)=
\dd A$.\qed\\

We shall use the following theorem, proven, for instance, in \cite{borsuk1969multidimensional}.\\

\begin{theorem}\label{Hyperplane intersection}\textsc{(Hyperplane intersection.)}
	The intersection of two hyperplanes in $\mathbb{R}^n$ with dimensions $k'$ and $k''$ is either the empty set or a hyperplane of dimension equal to at least $k'+k''-n$.\\
\end{theorem}

\begin{lm}
	\label{point dodging}
	\textsc{(Point dodging.)}
	Let $a$ and $b$ be two distinct points in $\mathcal{R}^2$ and $A$ be a finite set of points in $\mathcal{R}^2$ not containing $a$ and $b$. Then there exists a simple curve in $\mathbb{R}^2$ with endpoints $a$ and $b$, which is incident with no point of $A$.
\end{lm}

\textit{Proof}. Let $\{U_x\mid x\in A\}$ be a family of mutually disjoint closed disks none of which contains $a$ or $b$ and for each $x\in A$ the centre of $U_x$ is $x$. Let $B\lh\{x\in A\mid x\in[a,b]\}$ be the set of those points of $A$ that lie on the line segment $[a,b]$. Let $x\in B$. Notice that $[a,b]\cap U_x$ is a diameter of $U_x$. Let $a_x$ and $b_x$ be its endpoints. Let $\gamma$ be the curve obtained from the line segment $[a,b]$ by substituting each such diameter $[a_x,b_x]$ with some arc of $U_x$ with endpoints $a_x$ and $b_x$. Evidently $\gamma$ is a curve with the desired property.\qed\\

\begin{lm}
	\label{point dodging in connected open sets}
	\textsc{(Point dodging in connected open sets.)}
	Let $E$ be a connected and open in $\mathbb{R}^2$ set, $a$ and $b$ be two distinct points in $E$ and $A$ be a finite set of points in $\mathcal{R}^2$ not containing $a$ and $b$. Then there exists a curve contained in $E$ with endpoints $a$ and $b$ which is not incident with any point in $A$.
\end{lm}

\textit{Proof}. We know that a connected open in $\mathbb{R}^2$ set is homeomorphic to $\mathbb{R}^2$. Let $\phi$ be such a homeomorphism. By the point dodging lemma \ref{point dodging}, let $\gamma$ be a curve in $\mathbb{R}^2$ with endpoints $\phi(a)$ and $\phi(b)$ which is not incident with any point of $\phi[A]=\{\phi(x)\mid x\in A\}$. Then evidently the curve $\tilde\gamma\lh\{\la r,\phi^{-1}(\gamma(r))\ra\mid r\in Dom(\gamma)=[0,1]\}$ is a curve with the desired property.\qed\\

\begin{lm}
	\label{Dodging}
	\textsc{(Dodging.)}
	Let $n\geq 2 $ and $A$ be a finite set of $(n-2)$-dimensional hyperplanes in $\mathbb{R}^n$. Let $E$ be an open in $\mathbb{R}^n$ set and $a$ and $b$ be two points in $E\setminus(\cup A)$. Then there exists a simple curve contained in $E$ which is not incident with any element of $A$.
\end{lm}

\textit{Proof}. Induction on $n$.\\

Base: $n=2$. This is the point dodging in connected open sets lemma \ref{point dodging in connected open sets}. \\

Induction hypothesis: Let the claim be true for dimensions $k$ such that $2\leq k<n$.\\

Induction step: Let $A$ be a finite set of $(n-2)$-dimensional hyperplanes in $\mathbb{R}^n$ and $a$ and $b$ be points in $\mathbb{R}^n$ such that $a\notin\cup A$ and $b\notin\cup A$. Let $L$ be the set of all $(n-1)$-dimensional hyperplanes in $\mathbb{R}^n$ containing (the straight line connecting) $a$ and $b$. Clearly $|L|\geq\aleph_0$.\\

Let $\alpha\in A$ and $\lambda\in L$. Consider what $\alpha\cap\lambda$ could be. By the hyperplane intersection theorem, $\alpha\cap\lambda$ is either empty or a hyperplane of dimension $n-2$ or a hyperplane of dimension $n-3$. Evidently, since $\alpha$ is $(n-2)$-dimensional, $\alpha\cap\lambda$ is a hyperplane of dimension $n-2$ iff $\alpha\subseteq\lambda$.\\

We will show that for each $\alpha\in A$ there is at most one $\lambda\in L$ such that $\alpha\subseteq\lambda$. Suppose the contrary. Let $\alpha$, $\lambda_1$ and $\lambda_2$ be such. Then $\alpha\subseteq\lambda_1\cap\lambda_2$. Since $\lambda_1\neq\lambda_2$, by hyperplane intersection theorem $\lambda_1\cap\lambda_2$ has dimension $n-2$. Thus $\alpha=\lambda_1\cap\lambda_2$. But $a\in\lambda_1\cap\lambda_2$, which contradicts $a\notin\cup A$.\\

Thus only finitely many elements of $L$ have $(n-2)$-dimensional intersection with some element of $A$. But $L$ is infinite, so let $\lambda$ be an element of $L$ such that $B=\{\alpha\cap\lambda\mid\alpha\in A\}\setminus\{\emptyset\}$ is a finite set of $(n-3)$-dimensional hyperplanes.\\

By the induction hypothesis, let $\gamma$ be a simple curve in $E\cap\lambda$ with endpoints $a$ and $b$ which is not incident with any element of $B$, i.e. such that $\Ra(\gamma)\cap(\cup B)=\emptyset$. Then $\gamma$ is a simple curve in $\mathbb{R}^n$ with endpoints $a$ and $b$ which is not incident with any element of $A$.\qed\\

\begin{lm}
	\label{Infinity}
	\textsc{(Infinity.)}
	Let $n\geq 2$, $A$ be a polytope in $\mathbb{R}^n$ and $E$ be a connected open set such that $E\cap\dd A\neq\emptyset$. Then $|E\cap\dd A|\geq\aleph_0$.
\end{lm}

\textit{Proof}. Suppose $|E\cap\dd A|<\aleph_0$. Then $E\cap\dd A$ is a finite set of isolated points. Let $x\in E\cap\dd A$. Since $A$ is regular closed, let $a\in E\cap\Int(A)$ and $b\in E\setminus A$. By the dodging lemma \ref{Dodging}, there exists a simple curve contained in $E$ with endpoints $a$ and $b$ which is not incident with any point of $E\cap\dd A$. Let $\gamma$ be such. Then $\Ra(\gamma)\cap\dd A=\emptyset$, which contradicts the crossing lemma \ref{crossing}. Thus indeed $|E\cap\dd A|\geq\aleph_0$. \qed\\

\begin{lm}
	\label{Distributivity2}
	\textsc{(Distributivity.)}
	Let $A$, $B$ and $D$ be polytopes in $\mathbb{R}^2$ and $SC(A,B\cup D)$. Then $SC(A,B)$ or $SC(A,D)$.
\end{lm}

\textit{Proof}.\\

Case 1 : $A\sqcap(B\cup D)\neq\emptyset$. I.e. $A$ and $B\cup D$ overlap. Since the overlap relation is a contact relation, it is distributive over the join $\cup$. Thus $A\sqcap B=\emptyset$ or $A\sqcap D=\emptyset$. Then, by the upward strength lemma \ref{upward strength}, $SC(A,B)$ or $SC(A,D)$.\\

Case 2 : $A\sqcap(B\cup D)=\emptyset$. \\

Let us designate $B\cup D$ by $G$. Then $\Cl(\Int(A\cap G))=\emptyset$, thus $\Int(A\cap G)=\Int(A)\cap\Int(G)=\emptyset$. 
Let $E$ be a witness to $SC(A,G)$.\\

Since a polytope in $\mathbb{R}^2$ is a finite union of finite regularised intersections of closed half-planes, let $A=\cup_i\sqcap_j\alpha_{ij}$, $B=\cup_i\sqcap_j\beta_{ij}$ and $D=\cup_i\sqcap_j\delta_{ij}$, where the various $\alpha_{ij}$, $\beta_{ij}$ and $\delta_{ij}$ are closed half-planes and the indices vary through some six finite index sets. Let $P^A\lh\{\alpha_{ij}\mid i,j\}$, $P^B\lh\{\beta_{ij}\mid i,j\}$, $P^D\lh\{\delta_{ij}\mid i,j\}$ and $P\lh P^A\cup P^B\cup P^D$. Let $Q^A\lh\{\partial\alpha_{ij}\mid i,j\}$, $Q^B\lh\{\partial\beta_{ij}\mid i,j\}$, $Q^D\lh\{\partial\delta_{ij}\mid i,j\}$ and $Q\lh Q^A\cup Q^B\cup Q^D$. Then $Q^A$, $Q^B$, $Q^D$ and $Q$ are finite sets of lines.\\

Let us point out that $\dd A 			~=~ 
\dd(\cup_i\sqcap_j\alpha_{ij})			~\subseteq~
\cup_i\dd(\sqcap_j\alpha_{ij})			~=~
\cup_i\dd(\Cl(\Int(\cap_j\alpha_{ij})))	~\subseteq~
\cup_i\dd(\Int(\cap_j\alpha_{ij}))		~\subseteq~
\cup_i\dd(\cap_j\alpha_{ij})			~\subseteq~
\cup_i\cup_j\dd\alpha_{ij}				~=~
\cup Q^A$ and analogically for $B$ and $D$.\\

In the first half of the remaining part of the proof, we will show that there exists a point in $E\cap\dd A$ which is incident with exactly one element of $Q$. In the second half we will construct a sufficiently small open disk with centre such a point and will show that it is a witness to $SC(A,B)$ or to $SC(A,D)$.\\

First, we shall prove that $(E\cap\dd A) \cup (E\cap\dd G)\neq\emptyset$. Suppose the contrary, i.e. $E\cap\dd A=\emptyset$ and $E\cap\dd G=\emptyset$. Since $E\subseteq A\cup G$, we obtain $E=E\cap(A\cup G)=(E\cap A)\cup (A\cap G)=(E\cap(\Int(A)\cup\dd A))\cup(E\cap(\Int(G)\cup\dd G))=(E\cap\Int(A))\cup(E\cap\Int(G))\cup(E\cap\dd A)\cup(E\cap\dd G)=(E\cap\Int(A))\cup(E\cap\Int(G))$. But since $\Int(A\cap G)=\emptyset$, we have that $E\cap\Int(A\cap G)=E\cap(\Int(A)\cap\Int(G))=(E\cap\Int(A))\cap(E\cap\Int(G))=\emptyset$. Thus we obtained that $E$ is the union of two open disjoint sets, i.e. that $E$ is not connected, which is a contradiction. Thus indeed $(E\cap\dd A )\cup( E\cap\dd G)\neq\emptyset$.\\

Now we shall prove that $E\cap\dd A=E\cap\dd G$. Let $x\in E\cap\dd G$. Since $G$ is regular closed, let $a\in E\setminus G$ and $b\in E\cap\Int(G)$. Since $E\subseteq A\cup G$, we have $a\in A$. Suppose $b\in A$. Let $U$ be an open neighbourhood of $b$ contained in $\Int(G)$. Since $A$ is regular closed, we have $U\cap\Int(A)\cap\Int(G)\neq\emptyset$, which contradicts $\Int(A)\cap\Int(G)=\emptyset$. Thus $b\notin A$. Then $a$ and $b$ are witnesses to the fact that $x\in E\cap\dd A$. But $x$ was an arbitrary element of $E\cap\dd G$, thus we conclude that $E\cap\dd G\subseteq E\cap\dd A$. Analogically we obtain that $E\cap\dd A\subseteq E\cap\dd G$. Thus indeed $E\cap\dd A=E\cap\dd G\neq\emptyset$.\\

We shall now prove that there exists $\mu\in Q^A$ such that $|E\cap\mu\cap\dd A|\geq\aleph_0$. 
Suppose the contrary, i.e. suppose $(\forall\mu\in Q^A)(|E\cap\mu\cap\partial A|<\aleph_0)$. 
We have that $|E\cap(\cup Q^A)\cap\partial A|=|\cup\{E\cap\mu\cap\partial A\mid \mu\in Q^A\}|\leq\Sigma_{\mu\in Q^A}|E\cap\mu\cap\dd A|$. But the last is a finite sum of natural numbers, thus is finite.
Thus $|E\cap(\cup Q^A)\cap\partial A|<\aleph_0$. But $\partial A\subseteq\cup Q^A$, thus $(E\cap(\cup Q^A)\cap\partial A)=(E\cap\partial A)$. Thus $|E\cap\partial A|<\aleph_0$ which contradicts the infinity lemma \ref{Infinity}. Thus there indeed exists $\mu\in Q^A$ such that $|E\cap\mu\cap\dd A|\geq\aleph_0$. Let $\dd\alpha$ be such.\\ 

Let $Q_{\alpha}\lh Q\setminus\{\partial\alpha,\partial\alpha^*\}=Q\setminus\{\partial\alpha\}$ and $A(\alpha)\lh E\cap\dd\alpha\cap\dd A$.\\

We shall prove that there exists a point of $A(\alpha)$ which belongs to no element of $Q$ other than $\dd\alpha$. I.e. that $(\exists y\in A(\alpha))(\neg\exists\mu\in Q_{\alpha})(y\in\mu)$, i.e. that $A(\alpha)\nsubseteq\cup Q_{\alpha}$. Suppose the contrary. I.e. suppose $(\forall y\in A(\alpha))(\exists\mu\in Q_{\alpha})(y\in\mu)$. Let $M$ be a choice function that provides witnesses to these existences, i.e. let $M:A(\alpha)\longrightarrow Q_{\alpha}$ such that $(\forall y\in A(\alpha))(M(y)\in Q_{\alpha} ~~\&~~ y\in M(y))$.

We shall prove that $M$ is injective. Let $y_1,y_2\in A(\alpha)$ and $y_1\neq y_2$. Suppose $ M(y_1)= M(y_2)\rightleftharpoons\mu$. Then $\mu$ is the unique straight line incident with both $y_1$ and $y_2$. But $y_1$ and $y_2$ are elements of $A(\alpha)=E\cap\dd\alpha\cap\dd A$, thus they both lie on the line $\partial\alpha$. Thus $\partial\alpha=\mu=M(y_1)=M(y_2)$. But $\mu\in Q_{\alpha}=Q\setminus\{\dd\alpha\}$, thus $\mu\neq\partial\alpha$, which is a contradiction. Thus $M$ is indeed injective.

But the injectivity of $M$ implies that $|A(\alpha)|\leq|Q_{\alpha}|$, which is a contradiction because $|A(\alpha)|\geq\aleph_0$ and $Q_{\alpha}$ is finite. Thus the assumption that $A(\alpha)\subseteq\cup Q_{\alpha}$ is not true. So let $x$ be such that $x\in A(\alpha)=E\cap\dd\alpha\cap\dd A$ and $(\forall\mu\in Q_{\alpha})(x\notin\mu)$. In other words, $x$ is a point of $E\cap\dd A$ which belongs to exactly one element of $Q$ -- the element $\dd\alpha$. \\

Let $\rho$ be the Euclidean distance in $\mathbb{R}^2$. Let $R\leftrightharpoons\{\rho(x,\mu)\mid \mu\in Q_{\alpha}\}$. Notice that since $x$ is not incident with any line in $Q_{\alpha}$, $R$ is a finite set of strictly positive numbers, thus has a nonzero minimum. Let $e\leftrightharpoons\rho(x,\dd E)$. Since $x\in E$ and $E$ is an open set, $e$ is also nonzero. Let $r\leftrightharpoons\frac{1}{2}\textit{min}(R\cup\{e\})$ and $U$ be the open disk with centre $x$ and radius $r$. Let $p\leftrightharpoons \dd\alpha\cap U$. Clearly $p$ is a diameter of $U$. Let $U_1$ and $U_2$ be the two open half-disks that $p$ divides $U$ into. Clearly $p$, $U_1$ and $U_2$ are disjoint and $U=p\cup U_1\cup U_2$. By the definition of $U$, we have that $(\cup Q_{\alpha})\cap U=\emptyset$. And then, since $p\subseteq\mu=\dd\alpha$, we have $(\cup Q)\cap U_1=(\cup Q)\cap U_2=\emptyset$. \\

We have that $x\in \partial A$, $U$ is an open neighbourhood of $x$ and $A$ is regular closed, so let $a\in U\cap\Int(A)$. Since $p\subseteq\partial A$ and $\partial A$ and $Int(A)$ are disjoint, $a\notin p$, thus $a\in U_1$ or $a\in U_2$. WLoG let $a\in U_1$.\\

Suppose $U_1\nsubseteq Int(A)$. Let $a'\in U_1$ and $a'\notin Int(A)$. Then by the crossing lemma \ref{crossing}, we have $[a,a']\cap\partial A\neq\emptyset$, where $[a,a']$ is the line segment with endpoints $a$ and $a'$. Let $a''\in[a,a']\cap\partial A$. Since $U_1$ is a half-disk, it is convex, and thus $[a,a']\subseteq U_1$, so $a''\in U_1$. Thus $U_1\cap \partial A\neq\emptyset$, contradicts $(\cup Q)\cap U_1=\emptyset$, because $\dd A\subseteq\cup Q^A\subseteq \cup Q$. Thus $U_1\subseteq Int(A)$.\\

Since $U$ is an open neighbourhood of $x$ and $x\in\partial A$, there exists a point in $U$ that is not an element of $A$. Let $b$ be such. Since $U_1\subseteq Int(A)\subseteq A$ and $p\subseteq\partial A\subseteq A$, we have that $b\in U_2$. Since $b\in U_2\subseteq U\subseteq E\subseteq A\cup G$ and $b\notin A$, we have that $b\in G$. But $G=B\cup D$, so $b\in B$ or $b\in D$. WLoG let $b\in B$. \\

We obtain that $U_2\subseteq\Int(B)$ analogically to the way we obtained that $U_1\subseteq\Int(A)$.\\

We already know that $U_1\subseteq A$, $U_2\subseteq B$ and $p\subseteq\partial A\subseteq A$. Thus, since $U=p\cup U_1\cup U_2$, we have that $U\subseteq A\cup B$. Moreover $a$ and $b$ are witnesses to $U\cap A\neq\emptyset$ and $U\cap B_2\neq\emptyset$ respectively. And  obviously $U$, being an open disk, is connected and open. Thus $U$ is a witness to $SC(A,B)$.\qed\\

Thus, in view of remark \ref{ContIndeed} and the corollary \ref{UpwardStrengthCorollary} to the upward strength lemma \ref{upward strength},  $SC^{\mathbb{R}^2}$ is indeed distributive over the join $\cup$ in $\PP^2$. We have obtained that $SC^{\mathbb{R}^2}$ satisfies all of the conditions for being a contact relation in $\PP^2$. Thus $\la\PP^2,SC^{\mathbb{R}^2}\ra$ is a contact algebra.\\

\subsection{Higher Dimensions}

Let us suppose that throughout this section $n$ is a fixed natural number greater than 2, $\mathbb{V}=\la V,\tau\ra$ is an $n$-dimensional Euclidean space and $\Int$, $\Cl$ and $\sqcap$ designate the interior, closure and regularised intersection operators in $\mathbb{V}$. We shall use subscripts to designate the corresponding operators in other topological spaces.\\

Recall that if $A=\sqcap B$ is a nonempty basic polytope for some finite set $B$ of closed half-spaces of an Euclidean space, we have $A=\sqcap B=\cap B$.\\

\begin{lm}
	\label{Division}
	\textsc{(Division.)}
	Let $A$ be a finite set of closed half-spaces of $\mathbb{V}$ and $x\in\Int(\cap A)$. Let $\beta$ be a closed half-space of $\mathbb{V}$ such that $x\in\dd\beta$. Let $A_1\lh A\cup\{\beta\}$ and $A_2\lh A\cup\{\beta^*\}$. Then $\sqcap A_1$ and $\sqcap A_2$ are nonempty and $x\in\dd(\sqcap A_1)$ and $x\in\dd(\sqcap A_2)$.
\end{lm}

\textit{Proof}. Let $\rho$ be the Euclidean distance in $\mathbb{V}$ and let $r\lh \frac{1}{2}min\{\rho(x,\dd\alpha)\mid\alpha\in A\}$. Since $A$ is finite and $x\notin\dd \alpha$ for any $\alpha\in A$, we have that $r$ is positive. Let $U$ be the open $n$-dimensional ball with centre $x$ and radius $r$. Evidently $\dd\beta$ divides $U$ into two (nonempty) half-balls, i.e. $U_1\lh U\cap\Int(\beta)$ and $U_2\lh\cap\Int(\beta^*)$ are open hlaf-balls such that $U_1\subseteq\Int(\cap A_1)\subseteq\sqcap A_1$ and $U_1\subseteq\sqcap A_2$. Evidently $x\in\dd(\sqcap A_1)$ and $x\in\dd(\sqcap A_2)$.\qed

\begin{cor}
	\label{DivisonCor}\rm
	If $A$ is a finite set of closed half-spaces of $\mathbb{V}$, and $\alpha$ is a half-space of $\mathbb{V}$, then $\dd\alpha\cap\Int(\cap A)\subseteq\dd(\cap(A\cup\{\alpha\}))$ and $\dd\alpha\cap\Int(\cap A)\subseteq\dd(\cap(A\cup\{\alpha^*\}))$.\\
\end{cor}

To prove that $SC^{\mathbb{R}^n}$ is distributive over $\cup$ for polytopes in $\mathbb{R}^n$, we shall use a representation of the boundaries of polytopes, which we shall describe in this section. We shall prove that the boundary $\dd A$ of a polytope $A$ in an $n$-dimensional Euclidean space can be represented as a union $(\cup S)\cup K$ where $S$ is a finite set of open $(n-1)$-dimensional sets and $K$ is a subset of a finite union of $(n-2)$-dimensional hyperplanes in $\mathbb{R}^n$.\\

Let $\upphi$ be a finite set of $(n-1)$-dimensional hyperplanes in $\mathbb{V}$. We shall call such a set \textit{a set of cuts in $\mathbb{V}$}. Let $\mu$ be a cut in $\mathbb{V}$. There exist exactly two half-spaces $\alpha$ and $\alpha^*$ of $\mathbb{V}$ such that $\mu=\dd\alpha=\dd\alpha^*$. We shall call $\alpha$ and $\alpha^*$ \textit{the $\mathbb{V}$-sides of $\mu$}.\\

By $\bar{\upphi}$ we shall designate the set of the $\mathbb{V}$-sides of the elements of $\upphi$. We shall refer to the elements of $\bar{\upphi}$ as \textit{$\upphi$-$\mathbb{V}$-sides}. Evidently $(\forall\alpha\in\bar{\upphi})(\alpha^*\in\bar{\upphi})$ and $|\bar{\upphi}|=2|\upphi|$.\\

Let $s$ be a nonempty set of $\upphi$-$\mathbb{V}$-sides. We shall say that $s$ is \textit{$\upphi$-admissible} iff $\sqcap s\neq\emptyset$. Notice that this implies that $(\forall \alpha\in s)(\alpha^*\notin s)$ and $\sqcap s=\cap s$. We shall designate by $\upphi_a$ the set of $\upphi$-admissible sets.\\ 

We shall call a set $s$ of $\upphi$-$\mathbb{V}$-sides a \textit{$\upphi$-alternative} iff $s$ is $\upphi$-admissible and $(\forall\alpha\in \bar{\upphi})(\alpha\in s ~\textmd{ or }~ \alpha^*\in s)$. Evidently a $\upphi$-alternative is a set of exactly $|\upphi|$ half-spaces of $\mathbb{V}$. We shall designate by $\upphi_A$ the set of $\upphi$-alternatives.\\

For each $\upphi$-admissible set $s$ we shall call $\cap s$ a $\upphi$-block. By $\upphi_b$ we shall designate the set $\cap[\upphi_a]=\{\cap s\mid s\in \upphi_a\}$ of $\upphi$-blocks.\\

For each $\upphi$-alternative $s$ we shall call $\cap s$ a $\upphi$-brick. By $\upphi_B$ we shall designate the set $\cap[\upphi_A]=\{\cap s\mid s\in\upphi_A\}$ of $\upphi$-bricks.\\

For each $\upphi$-block $s$ we shall call $\Int(s)$ a $\upphi$-core. By $\upphi_C$ we shall designate the set $\Int[\upphi_B]=\{\Int(s)\mid s\in\upphi_B\}$ of $\upphi$-bricks. Notice that each $\upphi$-core is nonempty.\\

It is easy to see that each $\upphi$-core is the interior of a unique $\upphi$-brick, which is the (regularised) intersection of a unique $\upphi$-alternative.\\

\begin{lm}
	\label{}
	%\textsc{(.)}
	All $\upphi$-cores are mutually disjoint.
\end{lm}

\textit{Proof}. Let $A$ and $B$ be $\upphi$-cores and $A'$ and $B'$ the $\upphi$-alternatives such that $A=\Int(\cap A')=\cap\Int[A']$ and $B=\Int(\cap B')=\cap\Int[B']$. Let $A\neq B$. Let $\alpha$ be a witness to this inequality. WLoG, let $\alpha\in A'$ and $\alpha\notin B'$. Then $\alpha^*\notin A'$ and $\alpha^*\in B'$. Then $A\subseteq\Int(\alpha)$ and $B\subseteq\Int(\alpha^*)$. Thus $A\cap B=\emptyset$.\qed\\

\begin{lm}
	\label{Building bricks}
	\textsc{(Building bricks.)}
	Each $\upphi$-block is the union of a unique set of $\upphi$-bricks.
\end{lm}

\textit{Proof}. Let $\cap A$ be a $\upphi$-block for some $\upphi$-admissible set $A$. Induction on $q=|\upphi|-|A|$.\\

Base: $q=|\upphi|-|A|=0$. Then $|A|=|\upphi|$, thus $A$ is a $\upphi$-alternative and, thus $\cap A$ is itself a $\upphi$-brick.\\

Induction hypothesis: Let the claim be true for any $\upphi$-admissible set $B$ such that $|\upphi|-|B|\leq q$.\\

Induction step: Let $A$ be a $\upphi$-admissible set such that $|\upphi|-|A|=q+1$. Then $|A|=|\upphi|-q-1$, thus $|A|<|\upphi|$. Then $\dd[A]\subsetneqq\upphi$. Let $\dd\alpha$ be a witness to this, i.e. $\dd\alpha\in\upphi$ and $\dd\alpha\notin\dd[A]$. Then $\alpha\notin A$ and $\alpha^*\notin A$.\\

Evidently $\cap A=\sqcap A=\Cl(\Int(\sqcap A))=\Cl(\Int((\sqcap A)\cap(\alpha\cup\alpha^*)))=(\sqcap A)\sqcap(\alpha\cup\alpha^*)= ((\sqcap A)\sqcap\alpha)\cup((\sqcap A)\sqcap\alpha^*)=(\sqcap(A\cup\{\alpha\}))\cup(\sqcap(A\cup\{\alpha^*\}))$.\\ % $=\cup\{~\cap(A\cup\{\alpha\})~,~\cap(A\cup\{\alpha^*\})~\}$.\\

Let us designate $A\cup\{\alpha\}$ and $A\cup\{\alpha^*\}$ by $B_1$ and $B_2$ respectively. By the choice of $\alpha$ we have that $|B_1|=|B_2|=|A|+1$, and thus $|\upphi|-|B_1|=|\upphi|-|B_2|=q$.\\ %Each of $B_1$ or $B_2$ is either empty or is $\upphi$-admissible, thus the induction hypothesis applies to them. Evidently this gives us a representation of $\sqcap A$ of the desired kind.\\

Case 1: Exactly one of $\sqcap B_1$ and $\sqcap B_2$ is empty. WLoG let $\sqcap B_1\neq\emptyset$. Then by the induction hypothesis $\sqcap B_1 =\cup b=\cup\{\sqcap b^1,...,\sqcap b^{t}\}$ for some $\upphi$-alternatives $ b_1^1,..., b_1^{t_1}$ and $\sqcap A=(\sqcap B_1)\cup\emptyset=\sqcap B_1=\cup\{\sqcap b^1,...,\sqcap b^{t}\}$.\\

Case 2: None of $\sqcap B_1$ and $\sqcap B_2$ is empty. Then by the induction hypothesis $\sqcap B_1 =\cup b_1=\cup\{\sqcap b_1^1,...,\sqcap b_1^{t_1}\}$ and $\sqcap B_2=\cup b_2=\cup\{\sqcap b_2^1,...,\sqcap b_2^{t_2}\}$ for some $\upphi$-alternatives $ b_1^1$,..., $b_1^{t_1}$, $b_2^1$,..., $b_2^{t_2}$. Then $\sqcap A=(\sqcap B_1)\cup(\sqcap B_2)=(\cup b_1)\cup(\cup b_2)=\cup(b_1\cup b_2)=\cup(\{\sqcap b_1^1,...,\sqcap b_1^{t_1}\}\cup\{\sqcap b_2^1,...,\sqcap b_2^{t_2}\})=\cup\{\sqcap b_1^1,...,\sqcap b_1^{t_1},\sqcap b_2^1,...,\sqcap b_2^{t_2}\}$. \qed\\

Evidently, for any finite set $\upphi$ of cuts in $\mathbb{V}$, we have $(\cup\upphi)\cap(\cup\upphi_C)=\emptyset$. Moreover, if $x\in V\setminus(\cup\upphi)$, then $x\in\upphi_C$. Thus we have $V=(\cup\upphi)\cup(\cup\upphi_C)$. So $V$ is the union of the disjoint sets $\cup\upphi$ and $\cup\upphi_C$.\\

Let $\mu$ be an arbitrary element of the set $\upphi$ of cuts in $\mathbb{V}$. Let $\upmu$ designate the topological space with universe $\mu$ and topology -- the induced by $\mathbb{V}$ topology on $\mu$. Then $\upmu$ is an $(n-1)$-dimensional Euclidean space.\\

Consider the intersections of the elements of $\upphi$ with $\mu$. Let $\nu\in\upphi\setminus\{\mu\}$. If $\mu$ and $\nu$ are not parallel, then $\mu\cap\nu$ is an $(n-2)$-dimensional hyperplanes in $\mathbb{V}$, thus is a $(dim(\mu)-1)$-dimensional hyperplane in $\upmu$. And if $\mu\parallel\nu$, then $\mu\cap\nu=\emptyset$. Let $\upphi^{\mu}$ designate the set $\{\mu\cap\nu\mid\nu\in\upphi~~\&~~\nu\nparallel\mu\}$. Clearly $\upphi^{\mu}$ is a set of cuts in $\upmu$.\\

Consider the intersections of the elements of $\bar{\upphi}$ with $\mu$. Let $\alpha\in\bar{\upphi}$. If $\mu\nparallel\dd\alpha$, then $\mu\cap\alpha$ is a closed (in $\upmu$) half-space of $\upmu$ with boundary (in $\upmu$) $\mu\cap\dd\alpha$. If $\mu\parallel\dd\alpha$, then either $\mu=\dd\alpha$, or $\mu$ is disjoint with one of the sides of $\dd\alpha$ and is a subset of the interior if the other. Evidently the set $\{\mu\cap\alpha\mid\alpha\in\bar{\upphi}~~\&~~\dd\alpha\nparallel\mu\}$ is the set of $\upphi^{\mu}$-$\upmu$-sides.\\

We shall designate by $\overline{\upphi^{\mu}}$, $\upphi^{\mu}_a$, $\upphi^{\mu}_A$, $\upphi^{\mu}_b$, $\upphi^{\mu}_B$ and $\upphi^{\mu}_C$ the sets of $\upphi^{\mu}$-$\upmu$-sides, $\upphi^{\mu}$-admissible sets, $\upphi^{\mu}$-alternatives, $\upphi^{\mu}$-blocks, $\upphi^{\mu}$-bricks and $\upphi^{\mu}$-cores respectively.\\

Let $s$ be a $\upphi^{\mu}$-core for some element $\mu$ of $\upphi$. Then we shall say that $s$ is a $\upphi$-sheet. We shall designate by $\upphi_S$ the set $\cup\{\upphi^{\mu}_C\mid\mu\in\upphi\}$ of all $\upphi$-sheets.\\

By $\upphi_L$ we shall designate the set $\cup\{\upphi^{\mu}\mid\mu\in\upphi\}$ of all $(n-2)$-dimensional hyperplanes in $\mathbb{V}$ which are intersections of elements of $\upphi$. We shall call them $\upphi$-intersections.\\

Let $\mu\in\upphi$ and $s\in\upphi^{\mu}_a$. By $\hat{s}$ we shall designate the set $\{\alpha\in\bar{\upphi}\mid\mu\cap\alpha\in s\}$. %We shall call $\hat{s}$ the \textit{bulge of $s$}. 
By $\check{s}$ we shall designate the set $\{\alpha\in\bar{\upphi}\mid\mu\subseteq\Int(\alpha)\}$ of those $\mathbb{V}$-sides of the parallel to $\mu$ elements of $\upphi$ which contain $\mu$ in their interiors. Finally, by $\dot{s}$ we shall designate $\hat{s}\cup\check{s}$.\\

Let $\mu\in\upphi$ and $s\in\upphi^{\mu}_A$. Let $s_1\lh\{\mu_1\}\cup\dot{s}$ and $s_2\lh\{\mu_2\}\cup\dot{s}$, where $\mu_1$ and $\mu_2$ are the $\mathbb{V}$-sides of $\mu$. Evidently $s_1$ and $s_2$ are $\upphi$-alternatives, thus $\cap s_1=\sqcap s_1$ and $\cap s_2=\sqcap s_2$ are $\upphi$-bricks. We shall call $\sqcap s_1$ and $\sqcap s_2$ \textit{the $s$-toasts}. We have $(\sqcap s_1)\cup(\sqcap s_2) = (\sqcap(\dot{s}\cup\{\mu_1\})) \cup (\sqcap(\dot{s}\cup\{\mu_2\})) = ((\sqcap\dot{s})\sqcap\mu_1) \cup ((\sqcap\dot{s})\sqcap\mu_2) = (\sqcap\dot{s}) \sqcap (\mu_1\cup\mu_2) = \sqcap\dot{s}$. By the division lemma \ref{Division} we immediately obtain the following\\

\begin{lm}
	\label{Boundary sheets}
	\textsc{(Boundary sheets.)}
	Any $\upphi$-sheet is a subset of the boundaries of its toasts.\\
\end{lm}

\begin{lm}
	\label{Containment}
	\textsc{(Containment.)}
	Any $\upphi$-sheet is disjoint with any $\upphi$-brick other than its toasts.
\end{lm}

\textit{Proof}. Let $\mu\in\upphi$, $s\in\upphi^{\mu}_C$ and $\cap s_1$ and $\cap s_2$ be the $s$-toasts. Let $x\in s$. Let $t$ be a $\upphi$-alternative other than $s_1$ and $s_2$. Let $\alpha$ be a witness to their inequality, i.e. let $\alpha\in s_1$, $\alpha\in s_2$ and $\alpha\notin t$. Then $\alpha^*\notin s_1$, $\alpha^*\notin s_2$ and $\alpha^*\in t$. Evidently $\alpha\neq\mu_1$ and $\alpha\neq\mu_2$. Thus $\alpha\in\hat{s}\cup\check{s}$.

Case 1: $\alpha\in\check{s}$. Then $s\subseteq\Int(\alpha)$ and $\cap t=\cap(t\cup\{\alpha^*\})\subseteq\alpha^*$.

Case 2: $\alpha\in\hat{s}$. We have $s=\Int_{\upmu}(\sqcap_{\upmu}\{\mu\cap\beta\mid\beta\in\hat{s}\})=\Int_{\upmu}(\cap\{\mu\cap\beta\mid\beta\in\hat{s}\})=\cap\Int_{\upmu}[\{\mu\cap\beta\mid\beta\in\hat{s}\}]$. Then $s\subseteq\Int_{\upmu}(\mu\cap\alpha)=\mu\cap\Int(\alpha)\subseteq\Int(\alpha)$. And again $\cap t\subseteq\alpha^*$.\qed\\

\begin{lm}
	\label{Entirety}
	\textsc{(Entirety.)}
	Let $A$ be a finite set of $\upphi$-bricks and $s$ be a $\upphi$-sheet. Then either $s\subseteq\dd A$ or $s\cap A=\emptyset$.
\end{lm}

\textit{Proof}. Let $B$ be the set of $\upphi$-alternatives such that $A=\cup(\cap[B])$. Let $\cap s_1$ and $\cap s_2$ be the $s$-toasts.\\

Case 1: $s_1\notin B$ and $s_2\notin B$. By the containment lemma \ref{Containment}, $s$ is disjoint with any element of $\cap[B]$, thus $s$ is disjoint  with $A$.\\

Case 2: $s_1\in B$ and $s_2\in B$. Then $\cap s_1\in\cap[B]$ and $\cap s_2\in\cap[B]$. Then $(\cap s_1)\cup(\cap s_2)\subseteq A$. But $(\cap s_1)\cup(\cap s_2)=\cap\dot{s}$, so $\dot{s}\subseteq A$ and thus $\Int(\cap\dot{s})\subseteq\Int(A)$. By the corollary \ref{DivisonCor} to the division lemma \ref{Division}, we have $s\subseteq\Int(\cap\dot{s})$. Then $s\subseteq\Int(A)$. Thus $s\cap\dd A=\emptyset$.\\

Case 3: $s_1\in B$ and $s_2\notin B$. Evidently any open neighbourhood of any point of $s$ contains interior points of $\cap s_1$ and of $\cap s_2$. Since all $\upphi$-cores are mutually disjoint, it contains point from $A$ and points exterior to $A$. Thus $s\subseteq \dd A$.\qed\\

Since the boundaries of unions and of intersections are subsets of the unions of the boundaries of the respective sets, we have that the boundary of a union of $\upphi$-bricks is a subset of $\cup\upphi$. Also, $\cup\upphi = \cup\{\mu\mid \mu\in\upphi\} = \cup\{(\cup\upphi^{\mu}_C)\cup(\cup\upphi^{\mu})\mid\mu\in\upphi\} = (\cup\{\upphi^{\mu}_C\mid\mu\in\upphi\})\cup(\cup\{\upphi^{\mu}\mid\mu\in\upphi\})=(\cup\upphi_S)\cup(\cup\upphi_L)$.\\

Let $B$ be a finite set of $\upphi$-bricks and $A=\cup B$. By the entirety lemma \ref{Entirety}, let $S\lh\{s\in\upphi_S\mid s\subseteq\dd A\}$ and $S'\lh \upphi_S\setminus S=\{s\in\upphi_S\mid s\cap A=\emptyset\}$. Then $\dd A\subseteq(\cup S)\cup(\cup\upphi_L)$ and $\cup S\subseteq\dd A$. Let $K\lh(\dd A)\setminus(\cup S)$. Then $K\subseteq\cup\upphi_L$ and $\dd A=(\cup S)\cup K$.\\

We have just obtained a representation of the boundary of an arbitrary finite union of $\upphi$-bricks as a finite union of $\upphi$-sheets plus some subset of the union of the $\upphi$-intersections. We shall call this representation \textit{the $\upphi$-representation of $\dd A$}.\\

\subsubsection*{Distributivity in Higher Dimensions}

The following result is known, for instance form \cite{borsuk1969multidimensional}\\

\begin{lm}
	\label{CoveringLemma}
	\textsc{(Covering.)}
	%A finite union of $(n-1)$-dimensional hyperplanes in an $n$-dimensional Euclidean space $\mathbb{V}$ has no nonempty open in $\mathbb{V}$ subsets.
	A finite union of $(n-1)$-dimensional hyperplanes in an $n$-dimensional Euclidean space $\mathbb{V}$ is not a superset of any nonempty open in $\mathbb{V}$ set.\\
\end{lm}

\begin{lm}
	\label{DistributivityN}
	\textsc{(Distributivity.)}
	Let $n>2$. Let $A$, $B$ and $D$ be polytopes in $\mathbb{R}^n$ such that $SC(A,B\cup D)$. Then $SC(A,B)$ or $SC(A,D)$
\end{lm}

\textit{Proof}. Let us designate $B\cup D$ by $G$. Let $E$ be a witness to $SC(A,G)$. If $A\sqcap G\neq\emptyset$, then the proof is trivial as in the two-dimensional case with overlapping. So let $A\sqcap G=\emptyset$. Analogically to the two-dimensional case without overlapping, we have that $E\cap\dd A=E\cap\dd G\neq\emptyset$. Let $a$ and $b$ be points in $E$ such that $a\in\Int(A)$ and $b\in\Int(G)$.\\

Let $A=\cup_i \sqcap_j\alpha_{ij}$ for some closed half-spaces $\alpha_{ij}$ of $\mathbb{R}^n$. Let $\upphi\lh\cup\{\dd\alpha_{ij}\mid i,j\}$ be the set of the boundaries of those half-spaces and $\bar{\upphi}$ be the set of all half-spaces of $\mathbb{R}^n$ the boundaries of which are elements of $\upphi$. Evidently $\upphi$ is a set of cuts in $\mathbb{R}^n$ and $\bar{\upphi}$ is the set of $\upphi$-$\mathbb{R}^n$-sides. Then $A$ is a finite union of $\upphi$-blocks (the blocks $\{\sqcap_j\alpha_{ij}\mid i\}$), so by the building bricks lemma \ref{Building bricks}, $A$ is a finite union of $\upphi$-bricks.\\

Let $\dd A=(\cup S)\cup K$ be the $\upphi$-representation of the boundary of $A$. Then $S$ is a finite set of subsets of $(n-1)$-dimensional hyperplanes in $\mathbb{R}^n$ which are open in the induced by $\mathbb{R}^n$ topology on them and $K$ is some subset of the finite union $\cup\upphi_L$ of $(n-2)$-hyperplanes in $\mathbb{R}^n$.\\

By the dodging lemma \ref{Dodging} there exists a curve contained in $E$ with endpoints $a$ and $b$ which does not intersect $\cup\upphi_L$. Let $\gamma$ be such. By the crossing lemma \ref{crossing}, $\Ra(\gamma)\cap\dd A\neq\emptyset$. Then $\Ra(\gamma)\cap(\cup S)\neq\emptyset$. Let $s\in S$ such that $\Ra(\gamma)\cap s\neq\emptyset$.\\

Let $\mu$ be the $(n-1)$-dimensional hyperplane containing $s$ and $\upmu$ be the topological space with carrier $\mu$ and topology the induced by $\mathbb{R}^n$ topology on $\mu$. Then $E\cap s$ is an open in $\upmu$ set.\\

Let $B=\cup_i\sqcap_j\beta_{ij}$ and $D=\cup_i\sqcap_j\delta_{ij}$ for some half-spaces $\beta_{ij}$, $\delta_{ij}$ of $\mathbb{R}^n$. Let $\chi\lh\{\dd\alpha_{ij}\mid i,j\}\cup\{\dd\beta_{ij}\mid i,j\}\cup\{\dd\delta_{ij}\mid i,j\}$. Let $X=\{\mu_1\cap\mu_2\mid \mu_1\in\chi~~\&~~\mu_2\in\chi~~\&~~\mu_1\neq\mu_2\}\setminus\{\emptyset\}$ be the set of intersections of the nonparallel boundaries of the half-spaces by which $A$, $B$ and $D$ are constructed. Then $X$ is a finite set of $(n-2)$-dimensional hyperplanes in $\mathbb{R}^n$.\\

By the covering lemma \ref{CoveringLemma} we have $(E\cap s)\nsubseteq\cup X$. Let $x$ be a witness to this. Then $x\in E$ and $\mu$ is the only element of $\chi$ to which $x$ belongs. \\

Now, analogically to the two-dimensional case, we obtain that the open ball with centre $x$ and radius $\frac{1}{2}min\{\rho(x,\nu)\mid\nu\in\chi\cup\{\dd E\}\setminus\{\mu\}\}$ is a witness to $SC(A,B)$ or $SC(A,D)$.\qed\\

Thus, in view of remark \ref{ContIndeed} and the corollary \ref{UpwardStrengthCorollary} to the upward strength lemma \ref{upward strength}, for any $n>0$, $SC^{\mathbb{R}^n}$ is a contact relation in $\mathcal{P}^n$. We shall call it \textit{strong contact}. We shall designate the contact algebra $\la \PP^n,SC^{\mathbb{R}^n}\ra$ by $PSC^n$ and shall call it the \textit{strong-contact algebra of polytopes in $\mathbb{R}^n$}.\\

\subsection{Connectedness}

We say that a contact algebra is \textit{connected} iff any element $a$ of its carrier other than the zero and the unit is in contact with its complement.\\

\begin{theorem}
	\label{Connectedness}
	\textsc{(Connectedness.)}
	The strong-contact algebras of polytopes are connected.
\end{theorem}

\textit{Proof}. Let $A$ be a polytope in \RRnm such that $A\neq\emptyset$ and $A^*\neq\emptyset$. Obviously $\Rn$ itself is a witness to $SC(A,A^*)$.\qed\\

\section{The Logic of the Strong Contact}

\subsection{A Formal System}

We shall describe a standard formal system $\mathfrak{F}$ for connected contact algebras. \\

Let the \textit{alphabet} of the language $\mathcal{L}$ of $\mathfrak{F}$ consist of: a countable set $\textit{Ind}$ of individual variables, the equality symbol $\equiv$, the symbols $\neg$ and $\vee$ for the logical operators negation and disjunction respectively, the unary and binary function symbols $-$ and $+$ respectively for the Boolean complement and join, and the binary predicate symbol $C$ for the contact relation.\\

The \textit{terms} in $\mathcal{L}$ are finite words defined recursively as follows: the individual variables are terms and if $a$ and $b$ are terms, then $-a$ and $a\cdot b$ are terms.\\

The \textit{formulas} in $\mathcal{L}$ are finite words defined recursively as follows: if $a$ and $b$ are terms, then $a\equiv b$ and $C(a,b)$ are formulas; if $\f$ and $\psi$ are formulas, then $\neg\f$ and $\f\vee\psi$ are formulas.\\

Let us introduce some abbreviations of terms and formulas.\\

If $\f$ and $\psi$ are formulas in $\mathcal{L}$:\\
let $\f\wedge\psi$ abbreviate $\neg((\neg\f)\vee(\neg\psi))$\\
let $\f\Rightarrow\psi$ abbreviate $(\neg\f)\vee\psi$\\
let $\f\Leftrightarrow\psi$ abbreviate $(\f\Rightarrow\psi)~\wedge~(\psi\Rightarrow\f)$\\

If $a$ and $b$ are terms in $\mathcal{L}$:\\
let $a\leq b$ abbreviate $a+b\equiv b$\\
let $a\cdot b$ abbreviate $-((-a)+(-b))$\\
let $0$ abbreviate $a\cdot(-a)$\\
let $1$ abbreviate $-0$\\
let $a\not\equiv b$ abbreviate $\neg(a\equiv b)$\\
let $\top$ abbreviate $a\equiv a$\\
let $\perp$ abbreviate $a\not\equiv a$\\

Let $\mathfrak{F}$ have only one rule of inference - \textit{modus ponens} (MP).\\

Let $\mathfrak{F}$ have the following axiom schemes:\\

(1) A complete set of axiom schemes for the classical propositional logic\\

(2) A set of axiom schemes for Boolean algebras\\

(3) A set of axiom schemes for contact relations: if $a$, $b$ and $c$ are terms of $\mathcal{L}$, then the following formulas are axioms of $\mathfrak{F}$:

\begin{center}
	\begin{tabular}{ c }
		$ \neg C(0,a) $\\
		$ C(a,b+c)\Leftrightarrow (C(a,b) ~\vee~ C(a,c))$\\
		$ C(a,b)\Rightarrow C(b,a)$\\
		$ a\not\equiv 0\Rightarrow C(a,a)$\\
	\end{tabular}
\end{center}

(4) The axiom scheme of connectedness: if $a$ is a term in $\mathcal{L}$, then the following is an axiom of $\mathfrak{F}$:
\begin{center}
	$a\not\equiv 0 ~\Rightarrow~ (a\not\equiv 1 ~\Rightarrow~ C(a,-a))$\\
\end{center}

\subsection{Semantics}

A \textit{structure for $\mathcal{L}$} consists of a nonempty set $A$, called the carrier of $\mathcal{A}$ or the universe of $\mathcal{A}$, functions $-':A\longrightarrow A$ and $+':A\times A\longrightarrow A$ and a binary relation $C'\subseteq A\times A$, called \textit{interpretations in $\mathcal{A}$ of $-$, $+$ and $C$} respectrively. The logical symbol $\equiv$ is interpreted as the equality.\\

Let $\mathcal{A}$ be a structure for $\mathcal{L}$ with carrier $A$ and interpretations of $-$, $+$ and $C$ respectively $-'$, $+'$ and $C'$. A \textit{valuation of $\mathcal{L}$ in} $\mathcal{A}$ is a function $v:\textit{Ind}\longrightarrow A$ extended to all terms and formulas in $\mathcal{L}$ by recursion on their construction in the following way: \\
If $a$ and $b$ are terms in $\mathcal{L}$, then let

$v(-a)=-'(v(a))$

$v(a+b)=v(a)+'v(b)$

$v(a\equiv b)=\mathbb{T}$ ~~iff~~ $v(a)=v(b)$ 

$v(C(a,b))=\mathbb{T}$ ~~iff~~ $C'(v(a),v(b))$ ;\\
if $\f$ and $\psi$ are formulas in $\mathcal{L}$, then let:

$v(\neg\f)=\mathbb{T}$ ~~iff~~ $v(\f)=\mathbb{F}$

$v(\f\vee\psi)=\mathbb{T}$ ~~iff~~ $v(\f)=\mathbb{T}$ or $v(\psi)=\mathbb{T}$ ,\\
where $\mathbb{T}$ and $\mathbb{F}$ are special sets chosen to designate truth and falsity.\\

Let $\mathcal{A}$ be a structure for $\mathcal{L}$, $-'$, $+'$ and $C'$ be the interpretations in $\mathcal{A}$ of $-$, $+$ and $C$ respectively, $v$ be a valuation of $\mathcal{L}$ in $\mathcal{A}$ and $\f$ be a formula in $\mathcal{L}$. Let the expression $\la \mathcal{A},v\ra\vDash\f$ abbreviate $v(\f)=\mathbb{T}$. We will read this as '$\f$ is true in $\mathcal{A}$ under $v$'. If for every valuation $v'$ of $\mathcal{L}$ in $\mathcal{A}$ we have $\la\mathcal{A},v'\ra\vDash\f$, then we say that $\f$ is true in $\mathcal{A}$, which we designate by $\mathcal{A}\vDash\f$.\\

A structure for $\mathcal{L}$ in which all axioms of $\mathfrak{F}$ are true is called a \textit{model of $\mathfrak{F}$}. The models of $\mathfrak{F}$ are, by the choice of axioms, the connected contact algebras. \\

\subsubsection*{Kripke semantics}

We shall pay special attention to the particular case of set-theoretic contact algebras, because the adjacency spaces that induce them have some important properties of Kripke frames. In fact they are often called Kripke frames.\\

Let $\mathcal{A}$ be the set-theoretic contact algebra, induced by the adjacency space $\mathcal{F}=\la W,R\ra$. Recall that this implies that $R$ is reflexive and symmetric. If $v$ is a valuation of $\mathcal{L}$ in $\mathcal{A}$, we also say that $v$ is a valuation of $\mathcal{L}$ in $\mathcal{F}$ and call $\la\F,v\ra$ a \textit{Kripke model}. We introduce the expressions $\la\F,v\ra\vDash\f$ and $\F\vDash\f$, which we read as '$\f$ is true in $\F$ under $v$' and '$\f$ is true in $\F$', as abbreviations for $\la\mathcal{A},v\ra\vDash\f$ and $\mathcal{A}\vDash\f$ respectively.

\begin{defi}

Let $\F=\la W,R\ra $ and $\F'=\la W',R'\ra$ be adjacency spaces and $f$ be a surjective function from $W$ onto $W'$. We say that $f$ is a $p$-morphism from $\F$ to $\F'$ if the following conditions are satisfied:

\begin{center}
	\begin{tabular}{ r l }
		\rm{(p1)} & $ (\forall x,y\in W)(xRy~\rightarrow~ f(x)R'f(y)) $\\
		\rm{(p2)} & $ (\forall x',y'\!\in W')(x'R'y'~\rightarrow~(\exists x,y\in W)(f(x)=x' ~\&~ f(y)=y' ~\&~ xRy)) $\\
	\end{tabular}
\end{center}

\end{defi}

If there exists a $p$-morphism from $\F$ to $\F'$, then $\F$ is said to be a $p$-morphic preimage of $\F'$ and $\F'$ -- to be a $p$-morphic image of $\F$. It is easy to see that a composition of $p$-morphisms is a $p$-morphism.\\

Let $\la\F,v\ra$ and $\la\F',v'\ra$ be Kripke models.
We say that $f$ is a $p$-morphism from $\la F,v\ra$ to $\la F',v'\ra$ iff $f$ is a $p$-morphism from $\F$ to $\F'$ and for every variable $p\in \textit{Ind}$ and every element $x$ of $W$ we have $x\in v(p)$ iff $f(x)\in v'(p)$. In such a case we shall say that $\la F,v\ra$ is a $p$-morphic preimage of $\la F',v'\ra$.
Known results are the following lemmas.\\

\begin{lm}
	\label{Pmorph1}
	\textsc{(p-morphism, first.)}
	Let $\la\F,v\ra$ and $\la F',v'\ra$ be Kripke models
	and $f$ be a $p$-morphism from $\la\F,v\ra$ to $\la\F',v'\ra$. Then for every formula $\f$ in $\mathcal{L}$, we have $\la\F,v\ra\vDash\f$ iff $\la\F',v'\ra\vDash\f$.\\
\end{lm}

\begin{lm}
	\label{Pmorph2}
	\textsc{(p-morphism, second.)}
	Let $\F$ and $\F'$ be adjacency spaces, $f$ be a $p$-morphism from $\F$ to $\F'$ and $v'$ be a valuation of $\mathcal{L}$ in $\F'$. Then there exists a valuation $v$ in $\F$ such that $\la\F,v\ra$ is a $p$-morphic preimage of $\la\F',v'\ra$.
\end{lm}

\begin{cor}
	\label{Pmorph2Cor}\rm
	If a formula $\f$ in $\mathcal{L}$ is not true in an adjacency space $\F'$, then $\f$ is not true in any $p$-morphic preimage $\F$ of $\F'$.\\
\end{cor}

We shall make crucial use of the following

\begin{theorem}
	\label{General completeness}
	\textsc{(Completeness, general)}
	Let $\f$ be a formula in $\mathcal{L}$. Then the following are equivalent:
	
	\rm{(1)} $\f$ is a theorem of $\mathfrak{F}$
	
	\rm{(2)} $\f$ is true in all connected adjacency spaces
	
	\rm{(3)} $\f$ is true in all finite connected adjacency spaces
\end{theorem}

This theorem is proved in the paper \cite{balbiani2007modal}, where the authors consider a formal system which is clearly equivalent to $\mathfrak{F}$.\\

\begin{defi}\rm
	We shall define some graph-theoretic notions for adjacency spaces. Let $\F=\la W,R\ra$ be an adjacency space. A $k$-sequence $\{x_i\}_{i<k}$ of cells of $\F$ such that $k>0$ and for each $i<k-1$, $x_iRx_{i+1}$ and $x_i\neq x_{i+1}$ is called a \textit{path in $\mathcal{F}$ (from $x_0$ to $x_{k-1}$)}. A \textit{simple path in $\mathcal{F}$} is a path in $\mathcal{F}$ which is an injection. A \textit{simple cycle in $\F$} is a simple path $\{x_i\}_{i<k}$ in $\F$ such that $k>2$ and $x_0Rx_{k-1}$. A \textit{cycle in $\F$} is a path $\{x_i\}_{i<k}$ in $\F$ such that $x_0Rx_{k-1}$ and which contains a subsequence which is a simple cycle.
\end{defi}

Two cells of $\F$ are called \textit{connected in $\F$} iff there exists a path in $\F$ from one of them to the other. An adjacency space is called \textit{connected} iff any two of its cells are connected in it.\\

\begin{lm}
	\label{Connectedness Spaces iff Algebras}
	\textsc{(Connectedness.)}
	A finite adjacency space is connected iff the induced by it set-theoretic contact algebra is connected.
\end{lm}

\textit{Proof}. Let $\F=\la W,R\ra$ be an adjacency space and $\mathcal{A}=\la \PP(W),W\setminus~, \cup\ra$ be the induced by $\F$ set-theoretic contact algebra.\\

Suppose $\mathcal{A}$ is connected. Let $x$ and $y$ be cells of $\F$. Suppose there is no path in $\F$ from $x$ to $y$. Let $R'(x)$ and $R'(y)$ be the sets of cells to which there are paths in $\F$ from $x$ and $y$ respectively and let $R(x)\lh R'(x)\cup\{x\}$ and $R(y)\lh R'(y)\cup\{y\}$. Obviously $x\in R(x)$, $x\notin R(y)$, $y\in R(y)$ and $y\notin R(x)$, thus neither of $R(x)$ and $R(y)$ is empty or equal to $W$. Clearly $R(x)\cap R(y)=\emptyset$. Then $R(y)\subseteq W\setminus R(x)$. Since $\mathcal{A}$ is connected, $C_R(R(x),W\setminus R(x))$, i.e. $(\exists u\in R(x))(\exists v\in R(y))uRv$ which is a contradiction.\\

Suppose $\mathcal{A}$ is not connected. Let $a$ be a nonempty subset of $W$ unequal to $W$, such that $\neg C_R(a,W\setminus a)$, i.e. $(\forall x\in a)(\forall y\in W\setminus a)x\overline{R}y$. Let $x$ and $y$ be arbitrary elements of $a$ and $W\setminus a$ respectively.
Suppose $\pi=(x,...,y)$ is a path in $\F$ from $x$ to $y$. We will show that there exists $i\in\textit{Dom}(\pi)-1=k-1$ such that $\pi(i)\in a$ and $\pi(i+1)\notin a$. Suppose the contrary, i.e. that for each $i<k-1$, either both $\pi(i)$ and $\pi(i+1)$ are in $a$ or both are in $W\setminus a$. Since $\pi(0)=x\in a$ we can obviously prove by induction that $y\in a$, which would be a contradiction. Thus there exists $i<k-1$ such that $\pi(i)\in a$ and $\pi(i+1)\in W\setminus a$. But since $\pi$ is a path in $\F$, this means that $\pi(i)R\pi(i+1)$, which contradicts $\neg C_R(a,W\setminus a)$.\qed\\

Let $\pi$ be a simple cycle in $\F$ and $a$ be an element of $\textit{Range}(\pi)$. Clearly there are exactly two elements $b_1$ and $b_2$ of $\textit{Range}(\pi)$ other than $a$ such that $aRb_1$ and $aRb_2$. We shall call them \textit{the adjacent to $a$ cells in $\pi$}.\\

Let $\F=\la W,R\ra$ be an adjacency space, $\pi$ be a cycle in $\F$ and $(a,b)$ be a subpath of $\pi$, i.e. $\pi=(u_1,...,u_i,a,b,v_1,...,v_j)$ for some cells $u_1$,..., $u_i$, $v_1$,..., $v_j$ of $\F$. By $\pi_{ab}$ and $\pi_{ba}$ we will designate the cycles in $\F$  $(a,u_i,...,u_1,v_j,...,v_1,b)$ and $(b,v_1,...,v_j,u_1,...,u_i,a)$ respectively. Clearly $\pi_{ab}$ and $\pi_{ba}$ are paths in $\F$ from $a$ to $b$ and from $b$ to $a$ respectively.\\

\subsection{Completeness}

\subsubsection{Untying}

We shall suppose that throughout this section a finite connected adjacency space $\mathcal{F}=\la W,R\ra$ is fixed, and we shall examine some of its properties.

\begin{defi}\rm
Let $\F$ have cycles and $\pi$ be a simple cycle in $\F$. Let $a$ appear in $\pi$ and $b$ be one of the two adjacent to $a$ cells in $\pi$. Let $a'\notin W$. Let
\begin{center}
	\begin{tabular}{ c l }
		$W'$ & $ \lh~ W\cup\{a'\} $\\
		$R'$ & $ \lh~ (R\setminus\{\la a,b\ra,\la b,a\ra\})\cup\{\la a',b\ra , \la b,a'\ra , \la a',a'\ra\} $
	\end{tabular}
\end{center}
We call $\G=\la W',R'\ra$ the \textit{obtained from $\F$ by breaking $\pi$ at $a$ next to $b$} adjacency space.
\end{defi}

Let $\G$ be obtained from  $\F$ by breaking $\pi$ at $a$ next to $b$. Let $\mu$ be a path in $\F$ from $x$ to $y$, i.e. $\mu=(x,u_1,...,u_i,y)$ for some cells $u_1$,...,$u_i$ of $\F$. By $\tilde{\mu}$ we shall designate the sequence obtained from $\mu$ by substituting all subpaths $(a,b)$ and $(b,a)$ of $\mu$ with $\pi_{ab}$ and $\pi_{ba}$ respectively. Clearly $\tilde{\mu}$ is a path in $\G$ from $x$ to $y$.\\

Let $\{\G_i\}_{i<\omega}$ be a sequence of adjacency spaces defined by the following recursion:

Base: $\G_0\lh\F$

Recursion step: If $\G_k$ is acyclic, let $\G_{k+1}\lh \G_k$. If $\G_k$ contains a cycle, choose a simple cycle $\pi$ in $\G_k$, choose an element $a$ of $\textit{Range}(\pi)$ and one of the two adjacent to $a$ cells in $\pi$, which we shall designate by $b$. Then let $\G_{k+1}$ be the adjacency space obtained from $\F$ by breaking $\pi$ at $a$ next to $b$.\\

We shall call such a sequence an \textit{untying of $\F$}. Clearly an untying is a sequence of finite adjacency spaces. We will prove some additional properties of untyings.\\

\begin{lm}
	\label{Untying 1}
	\textsc{(Untying, first.)}
	Let $\{\G_i\}_{i<\omega}$ be an untying of $\F$. Then, for any $k<\omega$, if $\G_k$ has a cycle, $\G_{k+1}$ has strictly less simple cycles than $\G_k$.
\end{lm}

\textit{Proof}. Let $k<\omega$, $\G_k$ have a cycle and $\G_{k+1}$ be obtained from $\G_k$ by breaking $\pi$ at $a$ next to $b$. Then $\pi$ is a simple cycle in $\G_k$ but not in $\G_{k+1}$. It remains to show that no new simple cycles have been added, i.e. that each simple cycle in $\G_{k+1}$ is a simple cycle in $\G_k$. Let $\mu$ be a simple cycle in $\G_{k+1}$. We will show that $\mu$ is a simple cycle in $\G_k$. Since $a'$ is adjacent to only one cell -- $b$, it cannot appear in any simple cycle. Thus $a'$ does not appear in $\mu$. Then it is obvious from the definition of $R_{k+1}$ that $\mu$ is a simple cycle in $\G_k$.  \qed

\begin{cor}
	\label{Untying 1 Cor}\rm
	The number of simple cycles in an untying is strictly decreasing until at some point an acyclic adjacency space is constructed. Then, by the construction, all consecutive adjacency spaces are equal to it. Thus any untying of a finite connected adjacency space converges. We shall call the limit of an untying of $\mathcal{F}$ \textit{an untied version of $\F$}. Thus an untied version of a finite connected adjacency space is a finite connected acyclic adjacency space.\\
\end{cor}

\begin{lm}
	\label{Untying 2}
	\textsc{(Untying, second.)}
	Let $\{G_i\}_{i<\omega}$ be an untying of $\F$. Then, for any $k<\omega$, $\G_k$ is connected.
\end{lm}

\textit{Proof}. Induction on $k$. Base: $\G_0=\F$ is connected.

Induction hypothesis: Let $\G_k$ be connected.

Induction step: If $\G_{k+1}=\G_k$ the claim is trivially true. Let $\G_{k+1}=\la W_{k+1},R_{k+1}\ra$ be obtained from $\G_k$ by breaking the simple cycle $\pi$ at $a$ next to $b$. Let $x$ and $y$ be elements of $W_{k+1}$. We will show that $x$ and $y$ are connected in $\G_{k+1}$

Case 1: None of $x$ and $y$ equals $a'$. Then $x$ and $y$ are both elements of $W_k$. Since $\G_k$ is connected, let $\mu$ be a path in $\G_k$ from $x$ to $y$. Then $\tilde{\mu}$ is a path in $\G_{k+1}$ from $x$ to $y$.

Case 2: One of $x$ and $y$ equals $a'$. WLoG let $x=a'$ Let $\mu$ be a path in $\G_k$ from $b$ to $y$ Then the concatenation $(a')\ast\tilde{\mu}$ of $(a')$ and $\tilde{\mu}$ is a path in $\G_{k+1}$ from $a'$ to $y$, i.e. from $x$ to $y$.\qed

\begin{cor}\rm
	\label{Untying 2 Cor}
	An untied version of a finite connected adjacency space is connected.\\
\end{cor}

\begin{lm}
	\label{Untying 3}
	\textsc{(Untying, third.)}
	Let $\{G_i\}_{i<\omega}$ be an untying of $\F$. Then, for any $k<\omega$, $\G_{k+1}$ is a $p$-morphic preimage of $\G_k$.
\end{lm}

\textit{Proof}. Let $k<\omega$. If $\G_k$ is acyclic, then $\G_{k+1}=\G_k$, the claim is true for trivial reasons, so let $\G_{k+1}$ is obtained from $\G_k$ by breaking $\pi$ at $a$ next to $b$. Let $f=Id_W\cup\{\la a',a\ra\}$, We will show that $f$ is a $p$-morphism from $\G_{k+1}$ to $\G_k$, i.e. that $f$ satisfies the following conditions:
\begin{center}
	\begin{tabular}{ c l }
		(p1) & $ (\forall x,y\in W_{k+1})(\la x,y\ra\in R_{k+1}~\rightarrow~\la f(x),f(y)\ra\in R_n)$\\
		(p2) & $ (\forall x,y\in W_k)(\la x,y\ra\in R_k ~\rightarrow~ $\\
		&$(\exists x',y'\in W_{k+1})(f(x')=x ~\&~ f(y')=y ~\&~ \la x',y'\ra\in R_{k+1})) $
	\end{tabular}
\end{center}
Clearly (p1) is satisfied. For (p2), if $x=a$ and $y=b$, or vice versa, then $a'$ and $b$ are witnesses to what we want to prove. For any other $x$ and $y$, $x'=x$ and $y'=y$ are such witnesses.\qed

\begin{cor}
	\label{Untying 3 Cor}
	An untied version of a finite connected adjacency space $\F$ is a $p$-morphic preimage of $\F$.\\
\end{cor}

\begin{theorem}
	\label{Untying}
	\textsc{(Untying.)}
	Every finite connected adjacency space is a $p$-morphic image of a finite connected acyclic adjacency space.
\end{theorem}

\textit{Proof}. Let $\G$ be an untied version of $\F$. By the corollaries \ref{Untying 1 Cor}, \ref{Untying 2 Cor} and \ref{Untying 3 Cor}, to the first, second and third untying lemmas \ref{Untying 1}, \ref{Untying 2} and \ref{Untying 3}, we immediately obtain that $\G$ is a finite connected reflexive and symmetric $p$-morphic preimage of $\F$.\qed\\

\subsubsection{Projection}

We shall suppose that throughout this section a finite connected acyclic adjacency space $\mathcal{F}=\la W,R\ra$ is fixed, and we shall examine some of its properties. Let also an arbitrary cell $\alpha$ of $\F$ be fixed. \\

Let $L'=\{L_i\}_{i<\omega}$ be the sequence defined by the following recursion:

Base: Let $L_0\lh\{\alpha\}$ contain only the element $\alpha$.

Recursion step: Let $L_{k+1}\lh\{x\in W\setminus\cup\{L_i\mid i\leq k\}\mid(\exists y\in L_k)xRy\}$ be the set of those elements of $W$ that do not appear in $L_i$ for any $i\leq k$ and which are adjacent to some element of $L_k$.\\

We call the nonempty elements of the sequence $L'$ \textit{$\alpha$-levels of $\mathcal{F}$}. \\

The connectedness of $\mathcal{F}$ ensures that each cell of $\F$ appears in some level. The very construction of $L'$ ensures that no cell appears in two distinct levels. Since $\mathcal{F}$ is finite, the $\omega$-sequence $L'$ has a finite initial segment of nonempty elements (levels), followed only by empty ones. Let $L$ be that initial segment. \\

We shall call $L$ \textit{the $\alpha$-hierarchy of levels of $\mathcal{F}$}. If $x$ is a cell of $\mathcal{F}$, by $l^{\alpha}(x)$ we will designate the unique natural number $k$ such that $x\in L_k$. 

\begin{defi}\rm
We call $\#$ an \textit{$\alpha$-numeration of $\mathcal{F}$}, if $\#:W\twoheadrightarrowtail |W|$ and for any elements $x$ and $y$ of $W$, $l^{\alpha}(x)<l^{\alpha}(y)$ implies $\#(x)<\#(y)$. If $\#$ is an $\alpha$-numeration of $\mathcal{F}$, we call the inverse function $\#^{-1}$ of $\#$ an \textit{$\alpha$-listing of $\mathcal{F}$}. We say that a function is \textit{a numeration of $\F$} it is an $x$-numeration of $\F$ for some cell $x$ of $\F$. Analogically for listings.
\end{defi}

\begin{lm}
	\label{Numerations}
	\textsc{(Numerations.)}
	Let $\#$ be an $\alpha$-numeration of $\mathcal{F}$. Then $(\forall x\in W)(x\neq\alpha~\rightarrow~(\exists!y\in W)(\#(y)<\#(x) ~~\&~~ xRy))$.
\end{lm}

\textit{Proof.} Induction on $\#(x)$. Base: $\#(x)=0$, thus $x=\alpha$, thus the implication is trivially true.\\

I.h.: Let the claim be true for all $x'$ such that $\#(x')<\#(x)$.\\

I.s.: Let $x\neq\alpha$. Then $x\in L_j$ for some $j\gneqq0$. By the construction of the hierarchy $L$ of levels, $(\exists y\in L_{j-1})(xRy)$. Let $y$ be such. Then $l^{\alpha}(y)< l^{\alpha}(x)$. Since $\#$ is an $\alpha$-numeration of $\F$, $\#(y)<\#(x)$. Thus $y$ is a witness to the existence. \\

Now suppose $y$ and $y'$ be two distinct such cells, i.e. let $y'$ also be such that $\#(y')<\#(x)$ and $xRy'$. Then, by the induction hypothesis, we can construct paths from $y_1$ and from $y_2$ to $\alpha$. Let $\pi\ast(\alpha)=(y_1,...,\alpha)$ and $\mu\ast(\alpha)=(y_2,...,\alpha)$ be such. Let $\mu'$ be the path $\mu$ in the reverse direction. Then evidently $\pi\ast(\alpha)\ast\mu'\ast(x)$ is a cycle in $\F$, which is a contradiction.\qed\\

\begin{lm}
	\label{Paths}
	\textsc{(Paths.)}
	Let $x$ be a cell of $\F$, other than $\alpha$. Then there exists a unique simple path $\pi$ in $\mathcal{F}$ from $\alpha$ to $x$. Moreover, if $\#$ is an $\alpha$-numeration of $\mathcal{F}$, no cell $y$ of $\pi$ is such that $\#(x)<\#(y)$.
\end{lm}

\textit{Proof.} Let $\#$ be an $\alpha$-numeration of $\F$. Let $\in W\setminus\{\alpha\}$. Induction on $\#(x)$. Let the claim be true for all cells $x'$ of $\F$ such that $x'\neq\alpha$ and $\#(x')<\#(x)$.\\

By the numerations lemma \ref{Numerations}, let $y$ be the unique cell of $\F$ such that $\#(y)<\#(x)$ and $xRy$. By the induction hypothesis let $\pi=(\alpha,...,y$ be the unique simple path from $\alpha$ to $y$. By the induction hypothesis we also have that for every cell $z$ in $\pi$, we have $\#(z)\leq\#(y)\lneqq\#(x)$. Then $\pi\ast(x)$ is evidently a simple path of the kind we need.\\

Now suppose $\pi\ast(x)$ is not unique. Let $\mu\ast(x)$ be another such simple path. Let $\mu'$ be the path $\mu$ in the reverse direction. Then evidently $\pi\ast(x)\ast\mu'$ is a cycle in $\F$, which is a contradiction.\qed\\

Let $\#$ be an $\alpha$-listing of $\mathcal{F}$ and let $w$ designate $|W|$. We will recursively define a sequence $\{J_i\}_{i<\omega}$ of sequences of cells of $\F$, called a \textit{$\#$-arrangement sequence of $\F$} such that for each $k$, $J_k$ is a sequence with domain the smaller of $2k+1$ and $2w+1$ and contains the elements of $\textit{Range}(\#\upharpoonright(k+1))=\{\#(0),...,\#(k)\}$, i.e.  $J_k:min\{2k+1,2w+1\}\twoheadrightarrow \textit{Range}(\#\upharpoonright(k+1))$.\\% and a sequence $\{S_k\}_k$ of symmetric binary relations in $W$.\\

Base: Let $J_0\lh(\alpha)=\{\la0,\alpha\ra\}$ be a sequence of only one element - $\alpha$. Clearly, $J_0:\{0\}\twoheadrightarrow\{\alpha\}$, and since 

\begin{center}
	\begin{tabular}{l}
		$\textit{Dom}(J_0)=\{0\}=1=min\{2.0+1,2w+1\} ~~ \textmd{and}$\\
		$\textit{Range}(J_0)=\{\alpha\}=\{\#(0)\}=\textit{Range}(\#\upharpoonright1)=\textit{Range}(\#\upharpoonright(0+1)) ~,~~ \textmd{indeed}$\\
		$J_0:min\{2.0+1,2w+1\}\twoheadrightarrow \textit{Range}(\#\upharpoonright(0+1))~.$
	\end{tabular}
\end{center}

Recursion hypothesis: Let $J_k$ be defined such that $J_k:(2k+1)\cap(2w+1)\twoheadrightarrow \textit{Range}(\#\upharpoonright(k+1))$.\\

Recursion step: Case 1: $k\geq w$. Then 

\begin{center}
	\begin{tabular}{l}
		$\textit{Dom}(J_k)=min\{2k+1,2w+1\}=2w+1 ~, \textmd{and} $\\
		$\textit{Range}(J_k)=\textit{Range}(\#\upharpoonright(w+1))=\{\#(0),\#(1),...,\#(w)\})=W ~.\textmd{ Then} $\\
		$J_k:(2w+1)\twoheadrightarrow W$. Then let $J_{k+1}\leftrightharpoons J_k$.
	\end{tabular}
\end{center}

Case 2: $k\leq w$. Then $min\{2(k+1)+1,2w+1\}=2(k+1)+1$. Let $b\leftrightharpoons\#(k+1)$. Let $a$ be the unique element of $W$ such that $\#(a)\lneqq\#(b)$ and $aRb$. Since $\#(a)\lneqq\#(b)$, $a\in\textit{Range}(\#\upharpoonright(k+1))$. By the induction hypothesis $\textit{Range}(J_k)=\textit{Range}(\#\upharpoonright(k+1))$, thus $a\in\textit{Range}(J_k)$. Choose $i$ such that $J_k(i)=a$.\\

Then define $J_{k+1}$ to be the sequence with length $lh(J_k)+2=(2k+1)+2=2(k+1)+1$ obtained from $J_k$ by substituting the chosen occurrence of $a$ with consecutive occurrences of $a$, $b$ and $a$ again. More explicitly, if 

\begin{center}
	\begin{tabular}{l}
		$J_k=(\alpha,u_1,...,u_{i-1},a,u_{i+1},...,u_j)$, then let \\
		$J_{k+1}\leftrightharpoons(\alpha,u_1,...,u_{i-1},a,b,a,u_{i+1},...,u_j)$.
	\end{tabular}
\end{center}

I.e. if 
\begin{center}
	\begin{tabular}{l}
		$J_k=\{\langle 0 , \alpha \rangle ,..., \langle i-1 , u_{i-1} \rangle , \langle i , a \rangle , \langle i+1 , u_{i+1} \rangle ,..., \langle j , u_j \rangle\}$, then let \\
		
		$J_{k+1}\leftrightharpoons\{\langle 0 , \alpha \rangle ,..., \langle i\!-\!1 , u_{i-1} \rangle , \langle i , a \rangle , \langle i\!+\!1 , b \rangle , \langle i\!+\!2 , a \rangle , \langle i\!+\!2 , u_{i+1} \rangle ,..., \langle j\!+\!2 , u_j \rangle\}$.
	\end{tabular}
\end{center}

Obviously, $J_{k+1}:(2(k+1)+1)\cap(2w+1)\twoheadrightarrow \textit{Range}(\#\upharpoonright((k+1)+1))$\\

Since $W$ is finite, the sequence $\{J_i\}_{i<\omega}$ converges. Let $J\leftrightharpoons \lim_{k\rightarrow\omega}J_k$. Then $J:(2w+1)\twoheadrightarrow W$. We call $J$ a $\#$-arrangement of $\mathcal{F}$. We call a surjection $J':(2w+1)\twoheadrightarrow W$ an \textit{arrangement of $\mathcal{F}$} iff it is a $\#$-arrangement of $\mathcal{F}$ for some listing $\#$ of $\mathcal{F}$.\\

\begin{lm}
	\label{Adjacency 1}
	\textsc{(Adjacency, first.)}
	Let $\#$ be an $\alpha$-listing of $\F$ and $\{J_i\}_{i<\omega}$ be a $\#$-arrangement sequence of $\F$. Then 
	\begin{center}
		$\forall k(\forall x,y\in \textit{Range}(J_k))(x\neq y\rightarrow(xRy\leftrightarrow\exists i(\{J_k(i),J_k(i+1)\}=\{x,y\})))~.$
	\end{center}
\end{lm}

\textit{Proof.} Induction on $k$. Base: $k=0$. The claim is trivially true because no two elements of $\textit{Range}(J_0)=\{\alpha\}$ are unequal.\\

I.h.: Let the claim be true for all $k'\leq k$.\\

I.s.: Case 1: $k\geq w$. Then $J_{k+1}=J_k$ and by the induction hypothesis the claim is true.\\

Case 2: $k<w$. Let $b\leftrightharpoons\#(k+1)$. Then $\textit{Range}(J_{k+1})=\textit{Range}(J_k)\cup\{b\}$. Let $a$ be the unique, according to the numerations lemma \ref{Numerations}, element of $\textit{Range}(J_k)=\textit{Range}(\#\upharpoonright(k+1))$ such that $aRb$. Let the chosen on the $(k+1)$'th recursion step occurrence of $a$ to be substituted with $(a,b,a)$ be on $j$'th place, i.e. let $j$ be such that $J_k(j)=a$ and $J_{k+1}(j)=J_{k+1}(j+2)=a$ and $J_{k+1}(j+1)=b$, and $\forall i(j<i\leq2w ~\rightarrow~ J_{k+1}(i+2)=J_k(i))$. Let $x,y\in\textit{Range}(J_{k+1})$. \\

Case 2.1: None of $x$ and $y$ equals $b$. Then $x,y\in\textit{Range}(J_k)$. \\

($\rightarrow$) : Suppose $xRy$. By the induction hypothesis let $i$ be such that $\{J_k(i),J_k(i+1)\}=\{x,y\}$. By the construction of $J_{k+1}$ it is clear that, if $i<j$, then $J_{k+1}(i)=J_k(i)$ and $J_{k+1}(i+1)=J_k(i+1)$ and thus $\{J_{k+1}(i),J_{k+1}(i+1)\}=\{J_k(i),J_k(i+1)\}=\{x,y\}$. Since none of $x$ and $y$ equals $b$, we have $i\neq j$ and $i\neq j+1$. If $i\geq j+2$, then by the construction of $J_{k+1}$ we have $J_{k+1}(i+2)=J_k(i)$ and $J_{k+1}(i+3)=J_k(i+1)$, thus $\{J_{k+1}(i+2),J_{k+1}(i+3)\}=\{J_k(i),J_k(i+1)\}=\{x,y\}$. Thus in all possible cases we have a witness $i$ to what we need.\\ 

($\leftarrow$) : The proof in this direction in this case is completely analogical to the proof in the other direction that has just been carried out.\\

Case 2.2: One of $x$ and $y$ equals $b$. WLoG let $y=b$. Then, since $x\neq y=b$ and $x\in\textit{Range}(J_{k+1})=\textit{Range}(J_k)\cup\{b\}$, $x\in\textit{Range}(J_k)$. \\

($\rightarrow$) : Let $xRy$, i.e. $xRb$. By the numerations lemma \ref{Numerations} and since $\textit{Range}(J_{k+1})=\textit{Range}(\#(k+2))$ we obtain that $x=a$. Obviously by the construction of $J_{k+1}$ we have $\{J_{k+1}(j),J_{k+1}(j+1)\}=\{a,b\}=\{x,y\}$.\\

($\leftarrow$) : Let $x\overline{R}y$, i.e. $x\overline{R}b$. Then $x\neq a$. Since $b$ has only one occurrence in $J_{k+1}$ and it is surrounded by two occurrences of $a$, obviously $\neg\exists i(\{J_{k+1}(i),J_{k+1}(i+1)\}=\{x,b\}=\{x,y\})$.  \qed \\

Let $\#$ be an $\alpha$-listing of $\mathcal{F}$ and $J$ be an $\#$-arrangement of $\mathcal{F}$. 

\begin{defi}
	\label{ProjDef}
	\textsc{(Projection.)}\rm\\
	
	Let $f'$ be the function with domain $W$ mapping each element $x$ of $W$ to the union of exactly those closed intervals $[k,k+1]=\{u\in\mathcal{R}^1\mid k\leq u\leq k+1\}$ such that $J(k)=x$. I.e. for each element $x$ of $W$, let $f'(x)=\cup\{[k,k+1]\mid J(k)=x\}$.\\

	Let $f$ be the function with domain $W$ such that, for each element $x$ of $W\setminus\{\alpha\}$, $f(x)=f'(x)$ and $f(\alpha)=f'(\alpha)\cup(\mathcal{R}^1\setminus(\cup \textit{Range}(f')))=(-\infty,0)\cup f'(\alpha)\cup(2w+2,+\infty)$. We shall call such a function the \textit{$J$-projection of $\mathcal{F}$ onto $\mathbb{R}^1$}. We call a function a \textit{projection of $\mathcal{F}$ onto $\mathbb{R}^1$} if it is the $J$-projection of $\mathcal{F}$ onto $\mathbb{R}^1$ for some arrangement $J$ of $\mathcal{F}$.\\

	Let $f_n$ be the function with domain $W$ such that for each element $x$ of $W$, $f_n(x)=f(x)\times\mathcal{R}^{n-1}$ be the cylindrification of $f(x)$ to $\mathcal{R}^{n-1}$. We shall call such a function the \textit{$J$-projection of $\mathcal{F}$ onto $\mathbb{R}^n$}. We call a function a \textit{projection of $\mathcal{F}$ onto $\mathbb{R}^n$} if it is the $J$-projection of $\mathcal{F}$ onto $\mathbb{R}^n$ for some arrangement $J$ of $\mathcal{F}$. 

\end{defi}

\begin{rem}
	\label{Interiors Remark}
	\textsc{(Interiors.) }\rm
	It is obvious by the definition of $f$, that for any integer $k$, the open interval $(k,k+1)$ has nonempty intersection with the image $f(x)$ of precisely one element $x$ of $W$ and, moreover, that it is its subset. Analogically for the cylinders $(k,k+1)\times\mathcal{R}^{n-1}$ and $f_n$.\\
\end{rem}

\begin{lm}
	\label{Adjacency 2}
	\textsc{(Adjacency, second.)}
	Let $f$ be a projection of $\mathcal{F}$ onto $\mathbb{R}^1$. Then for any cells $x$ and $y$ of $\mathcal{F}$, $xRy$ iff $SC(f(x),f(y))$.
\end{lm}

\textit{Proof}. Let $J$ be an arrangement of $\F$ such that $f$ is a $J$-projection of $\F$ onto $\mathbb{R}^1$. Evidently if $x=y$ then we have both $xRy$ and $SC(f(x),f(y))$. So suppose $x\neq y$.\\

Let $xRy$. By the first adjacency lemma \ref{Adjacency 1}, WLoG let $i$ be such that $J(i-1)=x$ and $J(i)=y$. Then, since $f$ is a $J$-projection, $[i-1,i]\subseteq f(x)$ and $[i,i+1]\subseteq f(y)$. Clearly $[i-\frac{1}{2},i+\frac{1}{2}]$ is a witness to $SC(f(x),f(y))$.\\

Let $SC(f(x),f(y))$. Then, by the upward strength lemma \ref{upward strength}, $f(x)\cap f(y)\neq\emptyset$. Let $u\in f(x)\cap f(y)$. We have $u\in f(x)=\cup\{[k,k+1]\mid J(k)=x\}$ and $u\in f(y)=\cup\{[k,k+1]\mid J(k)=y\}$. Let $u\in[k_x,k_x+1]\subseteq f(x)$ and $u\in[k_y,k_y+1]\subseteq f(y)$. Since $x\neq y$ and $J$ is a function, we have $k_x\neq k_y$. WLoG let $k_x<k_y$. Since $[k_x,k_x+1]\cap[k_y,k_y+1]\neq\emptyset$, we conclude that $k_x+1=k_y=u$. Then, by the first adjacency lemma \ref{Adjacency 1}, $xRy$.\qed

\begin{cor}
	\label{Adjacency 2 Cor}\rm
	It is clear that this result holds for projection onto $\mathbb{R}^n$ as well. The witnesses there can be taken to be cylinders $\big[u-\frac{1}{2} , u+\frac{1}{2}\big]\times\mathcal{R}^{n-1}$, or the open balls with centre $u$ and radius $\frac{1}{2}$, instead of the open intervals $\big[u-\frac{1}{2},u+\frac{1}{2}\big]$.\\
\end{cor}

\begin{defi}
	\label{NpolAdjSp}\rm
	We shall call an adjacency space the carrier of which is the range of a projection onto $\mathbb{R}^n$ of a finite connected acyclic adjacency space and the adjacency relation of which is $SC^{\mathbb{R}^n}$ an \textit{$n$-polytope adjacency space}.\\
\end{defi}

\begin{theorem}
	\label{Projection}
	\textsc{(Projection.)}
	Every finite connected acyclic adjacency space is isomorphic to an $n$-polytope adjacency space.
\end{theorem}

\textit{Proof}. Let $\mathcal{F}=\la W,R\ra$ be a finite connected acyclic adjacency space and $f$ be a projection of $\mathcal{F}$ onto $\mathbb{R}^n$. Then $f$ is an injection of $W$ into $H$. By the corollary \ref{Adjacency 2 Cor} to the second adjacency lemma \ref{Adjacency 2}, for any elements $x$ and $y$ of $W$, $xRy$ iff $SC^{\mathbb{R}^n}(f(x),f(y))$. Thus $\mathcal{F}$ is isomorphic to the $n$-polytope adjacency space $\la\textit{Range}(f),SC^{\mathbb{R}^n}\ra$.\qed\\

\subsubsection{Merging}

Let $\F=\la W,R\ra$ be an $n$-polytope adjacency space.
Let $\mathcal{A}=\la B(W),C_R\ra$ be the induced by $\F$ contact algebra. We want to construct an isomorphic to $\mathcal{A}$ strong-contact algebra of polytopes in $\mathbb{R}^n$. We will show that the set-theoretic union $\cup$ maps $\mathcal{A}$ to such an algebra.
Let us designate the image $\cup[\PP(W)]$ of $\PP(W)$ under $\cup$ by $B$.  \\

\begin{lm}
	\label{Bijectivity}
	\textsc{(Bijectivity.)}
	$\cup$ is bijective from $\PP(W)$ to $B$.
\end{lm}

\textit{Proof}. $B$ is defined such that the surjectivity is obvious, thus we only have to show that it is injective. Let $a$ and $b$ be unequal subsets of $W$. Let $x$ be a witness to this inequality. WLoG let $x\in a$ and $x\notin b$. By the remark on interiors \ref{Interiors Remark}, let $k$ be such that $(k,k+1)\subseteq f(x)$ and $(k,k+1)\cap (\mathcal{R}^n\setminus f(x))=\emptyset$. Clearly $k+\frac{1}{2}\in\cup a$ and $k+\frac{1}{2}\notin\cup b$, thus $\cup a\neq\cup b$.\qed\\

\begin{lm}
	\label{Complement}
	\textsc{(Complement.)}
	Let $a$ be a subset of $W$. Then $\cup(W\setminus a)=(\cup a)^*$.
\end{lm}

\textit{Proof}. Let us designate $W\setminus a$ by $b$. Since $W$ is finite, $a$ and $b$ are finite. Let $a=\{x_1,...,x_k\}$ and $b=\{y_1,...,y_m\}$. Then $\cup a=x_1\cup...\cup x_k$ and $\cup b=y_1\cup...\cup y_m$ are polytopes. Obviously $\partial(\cup a)=\partial(\cup b)$. It is clear from the definition of a projection that $\Int(\cup a)$, $\partial(\cup a)$ and $\Int(\cup b)$ are disjoint and their union is $\mathcal{R}^n$.
Then $\cup(W\setminus a)=\cup b=\dd(\cup b)\cup \Int(\cup b)=\mathcal{R}^n\setminus\Int(\cup a)=\Cl(\mathcal{R}^n\setminus (\cup a))=(\cup a)^*$.\qed\\

\begin{lm}
	\label{Contact}
	\textsc{(Contact.)}
	For any subsets $a$ and $b$ of $W$, $C_R(a,b)$ iff $SC(\cup a,\cup b)$.
\end{lm}

\textit{Proof}. Let $a$ and $b$ be elements of $\PP(W)$. Then $a$ and $b$ are finite sets of polytopes, thus $\cup a$ and $\cup b$ are polytopes. \\

Suppose $C_R(a,b)$, i.e. $(\exists x\in a)(\exists y\in b)xRy$. Let $x$ and $y$ be witnesses to this, i.e. $x\in a$, $y\in b$ and $xRy$, i.e. $SC(x,y)$. Then $x\subseteq\cup a$ and $y\subseteq\cup b$ and by the monotony of $SC$ with respect to $\subseteq$ we obtain $SC(\cup a,\cup b)$. \\

Now suppose $SC(\cup a,\cup b)$. Let $a=\{x_1,...,x_p\}$ and $b=\{y_1,...,y_q\}$. By the distributivity of the strong contact over $\cup$, we obtain $SC(x_1,y_1)$ ~or~ ... ~or~ $SC(x_1,y_q)$ ~or~....~or~ $SC(x_p,y_1)$ ~or~...~or~ $SC(x_p,y_q)$. Let $SC(x_i,y_j)$ for some $i<p$ and $j<q$. Then $x_i$ and $y_j$ are witnesses to $C_R(a,b)$.\qed\\

\begin{theorem}
	\label{Merging}
	\textsc{(Merging.)}
	Every finite contact algebra induced by an $n$-polytope adjacency space is isomorphic to a subalgebra of the strong-contact algebra of polytopes in $\mathbb{R}^n$.
\end{theorem}

\textit{Proof}. Let $\F=\la W,R\ra$ be an $n$-polytope adjacency space and $\mathcal{A}$ be the contact algebra $\la\la \PP(W),W\setminus~,\cup\ra , C_R\ra$ induced by it. 
Trivially, for any sets $A$ and $B$ we have $\cup(A\cup B)=(\cup A)\cup(\cup B)$. By this, the bijectivity lemma \ref{Bijectivity}, the complement lemma \ref{Complement} and the contact lemma \ref{Contact}, we obtain that $\cup$ is an isomorphism from $\mathcal{A}=\la\la \PP(W),W\setminus~,\cup\ra,C_R\ra$ to $\la\la\cup[\PP(W)],*,\cup\ra,SC\ra$.\qed\\

\subsubsection{Completeness}

\begin{lm}
	\label{Subalgebra}
	\textsc{(Subalgebra.)}
	Let $\mathcal{A}$ and $\mathcal{B}$ be connected contact algebras and $\mathcal{A}$ be a subalgebra of $\mathcal{B}$. Let $\f$ be a formula in $\mathcal{L}$. Then, if $\f$ is not true in $\mathcal{A}$, then $\f$ is not true in $\mathcal{B}$.
\end{lm}

\textit{Proof}. Let $v$ be a witness that $\f$ is not true in $\mathcal{A}$, i.e. let $v$ be a valuation of $\mathcal{L}$ in $\mathcal{A}$ such that $\la\mathcal{A},v\ra\nvDash\f$. Then $v$ is also a valuation of $\mathcal{L}$ in $\mathcal{B}$. It is obvious that by induction on the construction of $\varphi$ we can obtain that $\la\mathcal{B},v\ra\nvDash\f$. Thus $\f$ is not true in $\mathcal{B}$.\qed\\

\begin{theorem}
	\label{Completeness}
	\textsc{(Completeness.)}
	Let $\f$ be a formula in $\mathcal{L}$ which is true in $PSC^n$. Then $\f$ is a theorem of $\mathfrak{F}$.
\end{theorem}

\textit{Proof}. Suppose $\f$ is not a theorem of $\mathfrak{F}$. By the general completeness theorem \ref{General completeness}, there exists a finite connected adjacency space in which $\f$ is not true. Let $\mathcal{F}$ be such. 
By the untying theorem \ref{Untying}, there exists a finite connected acyclic adjacency space which is a $p$-morphic preimage of $\mathcal{F}$. Let $\mathcal{G}$ be such. By the corollary \ref{Pmorph2Cor} to the second $p$-morphism lemma \ref{Pmorph2}, $\f$ is not true in $\mathcal{G}$. 
By the projection theorem \ref{Projection}, there exists an isomorphic to $\mathcal{G}$ $n$-polytope adjacency space. Let $\mathcal{H}$ be such. Then $\f$ is not true in $\mathcal{H}$. 
Let $\mathcal{A}$ be the induced by $\mathcal{H}$ set-theoretic contact algebra. Then $\f$ is not true in $\mathcal{A}$.
By the merging theorem \ref{Merging}, there exists an isomorphic to $\mathcal{A}$ subalgebra of the strong-contact algebra $PSC^n$ of polytopes in $\mathbb{R}^n$. Let $\mathcal{B}$ be such. Then $\f$ is not true in $\mathcal{B}$. Then, by the subalgebra lemma \ref{Subalgebra}, $\f$ is not true in $PSC^n$. \qed\\

\subsubsection{Standard Topological Contact}

With minor additional observations the same construction can be used to prove analogical completeness theorems for the standard topological contact $C$ (non-emptiness of the set-theoretic intersection) for the polytopes and for the regular closed in $\mathbb{R}^n$ sets.\\

By the definitions \ref{ProjDef} and \ref{NpolAdjSp} of projection and $n$-polytope adjacency space, it is evident that in $n$-polytope adjacency spaces the strong contact coincides with the standard topological contact $C$. I.e. if $\la W,R\ra$ is an $n$-polytope adjacency space for some $n$ and $x,y\in W$, then $SC(x,y)$ iff $x\cap y\neq\emptyset$. \\ 

By this we can immediately obtain a version of the contact lemma \ref{Contact} with the topological contact, namely: 

\begin{lm}
If $\la W,R\ra$ is an $n$-polytope adjacency space, then for any subsets $a$ and $b$ of $W$ we have $C_R(a,b)$ iff $C(\cup a,\cup b)$. 
\end{lm}

Indeed, $C_R(a,b)$ is equivalent to $(\exists x\in a)(\exists y\in b)SC(x,y)$ by the definition of $C_R$ which is equivalent to $(\exists x\in a)(\exists y\in b)(x\cap y\neq\emptyset)$, which means precisely $(\cup a)\cap(\cup b)\neq\emptyset$, i.e. $C(\cup a,\cup b)$. \\

Having this, we can immediately obtain a version of the merging theorem \ref{Merging} with the topological contact, namely: 

\begin{theorem}
	Every finite contact algebra induced by an $n$-polytope adjacency space is isomorphic to a subalgebra of the contact algebra $\la \PP^n,C\ra$ of polytopes in $\mathbb{R}^n$ with the topological contact. 
\end{theorem}

This allows us to obtain a version of the completeness theorem \ref{Completeness} with the standard topological contact. Namely: 

\begin{theorem}
	If a formula $\f$ in $\mathcal{L}$ is true in $\la \PP^n, C\ra$, then it is a theorem of $\mathfrak{F}$. 
\end{theorem}

Finally, since $\la \PP^n, C\ra$ is a subalgebra of $\la\mathcal{RC}(\mathbb{R}^n),C\ra$, by the subalgebra lemma \ref{Subalgebra} we immediately obtain the following

\begin{theorem}
	If a formula $\f$ in $\mathcal{L}$ is true in $\la \mathcal{RC}(\mathbb{R}^n), C\ra$, then it is a theorem of $\mathfrak{F}$.\\
\end{theorem}

\section{Conclusion}

We have defined a contact relation between polytopes, which is strictly stronger than the standard topological contact and strictly weaker than the overlap relation. We have proved that the universal fragments $L(\{PSC^n\})$ of the logics of the resultant strong-contact polytope algebras for arbitrary dimensions all coincide with the set $T(\mathfrak{F})$ of theorems of the standard quantifier-free formal system $\mathfrak{F}$ for connected contact algebras.\\

Moreover, we have that $T(\mathfrak{F})$ also coincides with the universal fragments $L(\{\la \PP^n,C\ra\})$ and $L(\{\la \mathcal{RC}(\mathbb{R}^n),C\ra\})$ of the logics respectively of the polytope algebras and algebras of regular closed in Euclidean spaces sets with the standard topological contact. From \cite{balbiani2007modal} we also know that $T(\mathfrak{F})$ coincides with the universal fragment $L(\{\la \mathcal{RC}(T),C\ra\mid T\in \mathcal{T}_{con}\})\rightleftharpoons L(\mathcal{T}_{con}^{~C})$ of the logic of the class of all algebras of regular closed sets in connected topological spaces, again with the topological contact. \\

In short, for any positive natural numbers $k$, $m$ and $n$, we have:
$$T(\mathfrak{F})=L(\{PSC^k\})=L(\{\la\PP^m,C\ra\})=L(\{\la\mathcal{RC}(\mathbb{R}^n),C\ra\})=L(\mathcal{T}_{con}^{~C})~.$$

In particular, we conclude that the quantifier-free language $\mathcal{L}(+,-,C)$ of $\mathfrak{F}$ cannot distinguish the strong contact from the topological contact for polytopes and cannot distinguish between dimensions of algebras of polytopes and of regular closed in Euclidean spaces sets.\\

\section*{List of some of the used abbreviations}

\begin{list}{~}{~}
	\item[] $\la a_0,...a_{n-1}\ra$ designates the ordered $n$-tuple of $a_0,...,a_{n-1}$ in the given order.
	\item[] $B(W)$ designates the set-theoretic Boolean algebra $\la \PP(W),W\setminus~,\cup\ra$ over the non-empty set $W$.
	\item[] $(a_0,a_1,...,a_{n-1})$ designates the $n$-sequence $\{\la0,a_0\ra,\la1,a_1\ra,...,\la n-1,a_{n-1}\ra\}$ of the sets $a_0$, $a_1$,..., $a_{n-1}$
	\item[] $|A|$ designates the cardinality of the set $A$
	\item[] $f[A]$ designates the image $\{f(x)\mid x\in A\}$ under the (class-)function $f$ of the subset $A$ of the domain $\textit{Dom}(f)$ of $f$
	\item[] $a\overline{R}b$ expresses that $a$ is not in the binary relation $R$ with $b$
	\item[] $\pi\ast \mu$ designates the concatenation of the sequences $\pi$ and $\mu$
\end{list}

%\nocite{*}
\bibliographystyle{plain}
\bibliography{arXiv1}

\end{document}